\documentclass[preprint,12pt]{elsarticle}

\usepackage{amssymb}
\usepackage{framed,multirow}
\usepackage{subfigure}
\usepackage{graphicx}%
\usepackage{multirow}%
\usepackage{amsmath,amssymb,amsfonts}%
\usepackage{amsthm}%
\usepackage{mathrsfs}%
\usepackage{appendix}%
\usepackage{xcolor}%
\usepackage{textcomp}%
\usepackage{manyfoot}%
\usepackage{booktabs}%
\usepackage{algorithm}%
\usepackage{algorithmicx}%
\usepackage{algpseudocode}%
\usepackage{listings}%
\usepackage{array}
\usepackage{amssymb}
\usepackage{latexsym}
\usepackage{hyperref}
\usepackage{marvosym}
\usepackage[thicklines]{cancel}

\theoremstyle{thmstyleone}%
\newtheorem{theorem}{Theorem}
\newtheorem{lemma}{Lemma}
\theoremstyle{thmstyletwo}%
\newtheorem{remark}{Remark}%

\theoremstyle{thmstylethree}%
\begin{document}

\begin{frontmatter}



\title{{\color{black}A Relaxed Lagrange Multiplier Approach for Phase Field Models}}

\author[label1]{Jinpeng Zhang} 

\affiliation[label1]{organization={Department  of  Mathematics,  University  of  Macau},
            city={Macao  SAR},
            country={China}}

\author[label2,label3]{Chaoyu Quan} 

\author[label2,label3]{Xiaoping Wang \corref{cor1}} 

\affiliation[label2]{organization={School of Science and Engineering, The Chinese University of Hong Kong (Shenzhen)},
addressline={Longgang},
city={Shenzhen},
state={Guangdong},
postcode={518172},
country={China}}  

\affiliation[label3]{organization={Shenzhen International Center for Industrial and Applied Mathematics, Shenzhen Research Institute of Big Data},
	city={Shenzhen},
	state={Guangdong},
    postcode={518172},
	country={China}}           
\cortext[cor1]{Corresponding author: wangxiaoping@cuhk.edu.cn}           

\begin{abstract}
	{\color{black} This paper introduces a novel relaxed Lagrange multiplier (RLM) method for designing efficient and energy-stable numerical schemes for {\color{black}phase field models}. The proposed approach reformulates the original model by introducing a time-dependent Lagrange multiplier $r(t)$, 
		whose evolution is governed by an ordinary differential equation featuring a relaxation parameter $\alpha > 0.$
		This key relaxation technique slows the dynamics of the multiplier, ensuring enhanced consistency between the modified and original systems after discretization, while simultaneously avoiding the need to solve nonlinear algebraic equations present in the original Lagrange multiplier (LM) method. We construct first- and second-order temporal discretizations based on the RLM approach and rigorously prove their unconditional stability for a modified energy. Furthermore, we establish the boundedness of the numerical Lagrange multiplier and demonstrate that it converges to $1$ as the relaxation parameter tends to zero and the modified energy converges to the original energy as ${\color{black}\tau}$ converges to zero.  Extensive numerical experiments confirm the theoretical findings, showcasing the method's second-order accuracy, unconditional energy stability, and significantly improved computational efficiency—requiring roughly half the cost of comparable scalar auxiliary variable (SAV) and LM methods. The RLM method effectively resolves the consistency issue of the SAV approach without requiring the nonlinear free energy to be bounded below, offering a robust and highly efficient alternative for simulating phase field models.}
	
\end{abstract}

\begin{keyword}
	Lagrange multiplier, Energy dissipation, Gradient flows, Phase field models
\end{keyword}

\end{frontmatter}
\section{Introduction}
Phase field models, conforming to the second law of thermodynamics, are widely used in various scientific and engineering fields, including materials science \cite{chen2002phase,chen2024phase,02elder,79leslie} and fluid dynamics \cite{98Anderson,12gao,17luo,1996gurtin}. Given the total free energy, they can be derived based on the dissipation mechanism and the variation of free energy. In this paper, we consider a free energy functional given by
\begin{align}
	E(\phi)=\int_{\Omega}\frac{1}{2}\phi\mathcal{L}\phi+F(\phi)d\boldsymbol{x}, \label{EO}
\end{align}
where $\Omega$ represents the spatial domain, $\mathcal{L}$ is a linear self-adjoint elliptic operator, and $F(\phi)$ is the nonlinear free energy density. The corresponding phase field model for the above free energy functional is then formulated as
\begin{align}
	\begin{split}
		\left \{
		\begin{array}{ll}
			\phi_t=-\mathcal{G}\mu,\\
			\mu=\frac{\delta E(\phi)}{\delta\phi}=\mathcal{L}\phi+F^{\prime}(\phi),
		\end{array}
		\right.
	\end{split}\label{GF}
\end{align}
with either periodic boundary conditions or homogeneous Neumann boundary conditions $\frac{\partial\phi}{\partial\mathbf{n}}|_{\partial\Omega}=\frac{\partial\mu}{\partial\mathbf{n}}|_{\partial\Omega}=0$, where $\mathbf{n}$ denotes the outward normal vector of $\partial\Omega$. Moreover, $\mathcal{G}$ is a semi-positive definite operator known as the mobility operator and $\mu$ is the chemical potential. By selecting different mobility operators $\mathcal{G}$ and free energy functionals $E(\phi)$, various phase field models can be derived, including the Allen-Cahn equation \cite{allen1979} and the Cahn-Hilliard equation \cite{cahn1958}.
With homogeneous Dirichlet, homogeneous Neumann or periodic boundary conditions, these models satisfy the following energy dissipation law:
\begin{align}
	\frac{d}{dt}E(\phi)=-\int_{\Omega}\frac{\delta E(\phi)}{\delta\phi}\mathcal{G}\frac{\delta E(\phi)}{\delta\phi}d\boldsymbol{x}\le 0.
\end{align}

Designing numerical schemes for phase field models presents several challenges, such as handling the nonlinear term $F(\phi)$ and preserving the energy dissipation law at the discrete level. In recent years, several popular numerical algorithms have been developed to address these challenges, such as the convex splitting approach \cite{eyre1998,baskaran2013,zhang2025two}, the stabilized linearly implicit approach \cite{zhu1999,shen2010,hao2020third}, the relaxation exponential Runge–Kutta method \cite{li2024relaxation}, the operator splitting approach \cite{li2022stability,li2022stability2,quan2026unconditional}, and the exponential time differencing (ETD) approach \cite{wang2016, fu2022energy, fu2024higher,hou2024energy,quan2025maximum}, along with the invariant energy quadratization (IEQ) approach \cite{yang2016,yang2017}, which can result in unconditionally energy stable, linear, second-order schemes for a wide range of gradient flow systems. However, the IEQ approach requires solving a coupled linear system with variable coefficients at each time step and assumes that the free energy density $F(\phi)$ is bounded from below.
As a remedy, the SAV approach was proposed in \cite{shen2018,shen2019}. It introduces scalar auxiliary variables to replace auxiliary function variables leading to a reformulated system that is mathematically equivalent to the original gradient flow system at the continuous level. The SAV approach retains the advantages of the IEQ approach and only assumes that the free energy $\int_{\Omega}F(\phi)d\boldsymbol{x}$ is bounded from below. It requires solving only two linear systems with constant coefficients at each time step. However, numerical errors can break the equivalence between the original system and the reformulated system. To address this problem, the relaxed scalar auxiliary variable (RSAV) method and its generalized version (R-GSAV) were recently proposed in \cite{jiang2022,zhang2022}. These approaches introduce a relaxation step into the original SAV framework to enhance the consistency between the original and reformulated systems, while retaining all the computational advantages of the SAV method.

{\color{black} Another promising approach is the Lagrange multiplier (LM) approach \cite{cheng2020}. Similar to the SAV method, it reformulates the original system into an equivalent modified model by introducing a Lagrange multiplier. A key advantage of this approach is its ability to rigorously preserve the stability of the original energy.  Furthermore, it is also shown that a unique solution to the nonlinear equation exists for sufficiently small time steps \cite{cheng2025unique}.}

{\color{black} In addition to the energy dissipation law, bound or positivity preservation is also an important issue in phase-field simulations. Representative methods include, but are not limited to, the cut-off approach~\cite{li2020arbitrarily,lu2013cutoff}, the implicit--explicit approach~\cite{tang2016implicit,liao2020energy}, the ETD method~\cite{du2019maximum,du2021maximum,quan2025maximum}, the convex splitting method~\cite{chen2019positivity}, the Lagrange multiplier approach~\cite{cheng2022new}, the integrating factor Runge--Kutta method~\cite{ju2021maximum,li2021stabilized}, the SAV-type method~\cite{huang2022bound,ju2022stabilized,ju2022generalized}, and the operator splitting method~\cite{li2022stability,li2022stability2}.}

	{\color{black} 
		In this paper, we propose a relaxed Lagrange multiplier (RLM) approach in which a Lagrange multiplier $r(t)$ is introduced to improve consistency between the modified and original systems. Instead of solving a non-linear algebraic equation for the Lagrange multiplier \cite{cheng2020}, our approach incorporates the derivative term $\frac{dr}{dt}$ and the constraint $\frac{d}{dt} \int_{\Omega} F(\phi)d\boldsymbol{x} = \int_{\Omega}f(\phi)\phi_{t}d\boldsymbol{x}$ via a relaxation parameter $\alpha$.  This provides more flexible control over the Lagrange multiplier $r$. The parameter $\alpha$ blackuces its rate of change, thus improving the consistency mentioned above.  This technique also avoids solving nonlinear equations at each time step and therefore improves computational efficiency and robustness. }
	{\color{black} As a trade-off, an additional term, \(\frac{r-1}{\alpha}\) (which is zero at the continuous PDE level) is added to the original energy.  We can then design a fully decoupled scheme that  preserves the energy decaying property for the modified energy. One can also show that the numerical Lagrange multiplier remains bounded, and that the additional term $\frac{r-1}{\alpha}$ converges to zero as $\tau \to 0$. }
	
	The advantages of the RLM method can be outlined as follows:
	\begin{itemize}
		\item It effectively resolves the inconsistency problem after discretization between the original system and the reformulated system.
		\item It does not require any assumption, unlike the original SAV approach that assumes the nonlinear free energy function $\int_{\Omega}F(\phi)d\boldsymbol{x}$ to be bounded from below.
		{\color{black} \item It does not require to solve any nonlinear equation, thereby avoiding issues related to the existence and uniqueness of numerical solutions.}
		\item It only requires to solve one linear system with constant coefficients, whereas the SAV approach and the LM approach requires to solve two linear systems. As a result, the computational cost of the RLM method is essentially half that of the LM approach and the SAV approach.
	\end{itemize}

	{\color{black}The remainder of this paper is organized as follows. In Section 2, we revisit the SAV and LM methods for phase field models. Subsequently, in Section 3,  we propose the {\color{black} RLM} method and rigorously prove its energy stability properties. In Section 4, we prove that the numerical Lagrange multiplier $r^n$ has lower and upper bounds, and then in Section 5, we discuss the influence of the stabilization parameter $\alpha$. Section 6 presents specific examples and numerical tests to verify the accuracy and effectiveness of the proposed RLM numerical schemes. Finally, we conclude with a brief summary.}

	{\color{black}\section{A brief review of the SAV and Lagrange multiplier methods}
		In this section, we will consider phase field models described by \eqref{GF} with periodic boundary condition and recall the SAV method \cite{shen2018} and the Lagrange multiplier method \cite{cheng2020}, to solve this system.
		
		\subsection{The original SAV method}
		
		Let $E_{0}(\phi)=\int_{\Omega}(F(\phi)-\frac{s}{2}\phi^{2})d\boldsymbol{x}$ and assume $E_{0}(\phi)+C_{0}>0$, where $C_{0}>0$ is a constant and $s\geq 0$ is a stabilization parameter. The key idea of the SAV method is to introduce a scalar auxiliary variable $\eta(t)=\sqrt{E_{0}(\phi)+C_{0}}$ and expand the original system \eqref{GF} into the following equivalent system, 
		\begin{align}
			&\phi_{t}=-\mathcal{G} \mu, \label{eq:SAV_c1}\\
			&\mu=\mathcal{L} \phi+s\phi+\frac{\eta(t)}{\sqrt{E_0(\phi)+C_0}}(f(\phi)-s\phi), \label{eq:SAV_c2}\\
			&\eta_t=\frac{1}{2 \sqrt{E_0(\phi)+C_0}} \int_{\Omega} (f(\phi)-s\phi) \phi_t d\boldsymbol{x} .\label{eq:SAV_c3}
		\end{align}
		For this reformulated system, the original energy can also be blackefined as
		\begin{align}
			E_{M}(\phi)=\int_{\Omega}\left(\frac{1}{2}\phi\mathcal{L}\phi+\frac{1}{2}s\phi^{2}\right)d\boldsymbol{x}+\eta^{2}-C_{0}. \label{eq:EM}
		\end{align}
		Instead of discretizing the original system ($\ref{GF}$), one can discretize the reformulated system ($\ref{eq:SAV_c1}$)-($\ref{eq:SAV_c3}$) to preserve the energy dissipation law.
		
		Consider the time domain $[0,T]$ and discretize it into uniform meshes $0=t_{0}<t_{1}<\cdots<t_{N}=T$, where $t_{n}=n{\color{black}\tau}$ and ${\color{black}\tau} = T/N$. Let $\psi^{n}$ represent the numerical approximation of the function $\psi(\cdot,t_{n})$.
		Using these notations, 
		a second-order SAV Crank-Nicolson (SAV-CN) semi-discrete scheme for ($\ref{eq:SAV_c1}$)-($\ref{eq:SAV_c3}$) can be obtained as follows:
		\begin{align}
			& \frac{\phi^{n+1}-\phi^n}{{{\color{black}\tau}}}=-\mathcal{G} \mu^{n+\frac{1}{2}}, \label{eq:SAV_d1}\\
			& \mu^{n+\frac{1}{2}}=\mathcal{L} \phi^{n+\frac{1}{2}}+s \phi^{n+\frac{1}{2}}+\frac{\eta^{n+\frac{1}{2}}}{\sqrt{E_0(\overline{\phi}^{n+\frac{1}{2}})+C_0}} \left(f(\overline{\phi}^{n+\frac{1}{2}})-s \overline{\phi}^{n+\frac{1}{2}}\right), \label{eq:SAV_d2}\\
			& \frac{\eta^{n+1}-\eta^n}{{{\color{black}\tau}}}=\int_{\Omega} \frac{\left(f(\overline{\phi}^{n+\frac{1}{2}})-s \overline{\phi}^{n+\frac{1}{2}}\right)}{2\sqrt{E_0(\overline{\phi}^{n+\frac{1}{2}})+C_0}} \frac{\phi^{n+1}-\phi^n}{{{\color{black}\tau}}} d\boldsymbol{x}, \label{eq:SAV_d3}
		\end{align}
		where $\phi^{n+\frac{1}{2}}=\frac{1}{2}(\phi^{n}+\phi^{n+1})$ and $\overline{\phi}^{n+\frac{1}{2}}=\frac{3}{2} \phi^n-\frac{1}{2} \phi^{n-1}$. This numerical scheme is also unconditionally energy stable, and its rigorous proof can be found in \cite{shen2018,shen2019}. \par
		The relaxed SAV (RSAV) method was proposed in \cite{jiang2022}. The RSAV method builds on the original SAV method by including\eqref{eq:SAV_d1}--\eqref{eq:SAV_d3} and introducing an additional correction step to improve the consistency between the original system and the reformulated system after discretization. For specific details, please refer to \cite{jiang2022}.
		
		\subsection{The Lagrange multiplier method (LM)}
		In the Lagrange multiplier method for phase field models \cite{cheng2020}, a constant Lagrange multiplier $q(t) \equiv 1$ is introduced to reformulate the original model \eqref{GF} as
		\begin{align}
			& \frac{\partial \phi}{\partial t}=-\mathcal{G} \mu, \label{eq:LM_c1}\\
			& \mu=\mathcal{L} \phi+q(t) F^{\prime}(\phi), \label{eq:LM_c2}\\
			& \frac{d}{d t} \int_{\Omega} F(\phi) d \boldsymbol{x}=q(t) \int_{\Omega} F^{\prime}(\phi) \phi_t d \boldsymbol{x}. \label{eq:LM_c3}
		\end{align}
		Instead of discretizing the original system $\eqref{GF}$, one can discretize the reformulated system \eqref{eq:LM_c1}--\eqref{eq:LM_c3}. A second-order LM-CN semi-discrete scheme for \eqref{eq:LM_c1}--\eqref{eq:LM_c3} can be obtained as follows:
		\begin{align}
			&\frac{\phi^{n+1}-\phi^n}{{\color{black}\tau}}=-\mathcal{G} \mu^{n+\frac{1}{2}}, \label{eq:LM_nu_c1}\\
			&\mu^{n+\frac{1}{2}}=\mathcal{L} \phi^{n+\frac{1}{2}}+F^{\prime}\left(\phi^{\star, n}\right) q^{n+\frac{1}{2}}, \label{eq:LM_nu_c2} \\
			&\left(F\left(\phi^{n+1}\right)-F\left(\phi^n\right), 1\right)=q^{n+\frac{1}{2}}\left(F^{\prime}\left(\phi^{\star, n}\right), \phi^{n+1}-\phi^n\right), \label{eq:LM_nu_c3}
		\end{align}
		where $q^{n+\frac{1}{2}}=\frac{1}{2}(q^{n}+q^{n+1})$ and $\phi^{\star, n}=\frac{3}{2} \phi^n-\frac{1}{2} \phi^{n-1}$. This numerical scheme is also unconditionally original energy stable, and its rigorous proof can be found in \cite{cheng2020}.
		\begin{remark}
			To maintain the unconditional stability of the original energy, it is necessary to solve the nonlinear algebraic equation \eqref{eq:LM_nu_c3}. As proven in \cite{cheng2025unique}, this equation has a unique solution when the time step is sufficiently small.
		\end{remark}

		\section{The relaxed Lagrange multiplier method  (RLM) for phase field models}
		
		{\color{black} We now introduce  a relaxed Lagrange multiplier method which allows us to decouple the solution of the Lagrange multiplier with the solution of the phase field model so as to avoid solving the nonlinear algebraic equation.
			We propose a relaxed version of (\ref{eq:LM_c3}) for the Lagrange multiplier $r(t)$ as follows:
			\begin{align}
				\frac{dr}{dt}=\alpha\left(-\frac{d}{dt}\int_{\Omega}F(\phi) d\boldsymbol{x}+r\int_{\Omega}f(\phi)\phi_{t}d\boldsymbol{x}\right),
			\end{align}
			where $\alpha >0 $ is an arbitrarily small relaxation parameter.
		} Taking into account the stabilization parameter $s$, we rewrite this evolution equation as follows:
		\begin{equation}
			\left \{
			\begin{aligned}
				\frac{dr}{dt}&=\alpha\left(-\frac{dE_{0}}{dt}+\int_{\Omega}(rf(\phi)-s\phi)\phi_{t}d\boldsymbol{x}\right),\\
				r(0)&=1,
			\end{aligned}
			\right.
			\label{CSAV_c4}
		\end{equation}
		where $E_{0}(\phi)=\int_{\Omega}(F(\phi)-\frac{s}{2}\phi^{2})d\boldsymbol{x}$. Intuitively, when $\alpha$ gets smaller,  the change of $r$ is slowed down. With this evolution equation of $r$, we can reformulate the system \eqref{GF} into the following equivalent form:}
	{ 
		\begin{align}
			& \phi_t=-\mathcal{G}\mu, \label{MCH1}\\
			& {\mu=\mathcal{L}\phi+s\phi+(rf(\phi)-s\phi),} \label{MCH2}\\
			& {\frac{dr}{dt}=\alpha\left(-\frac{dE_{0}(\phi)}{dt}+\int_{\Omega}(rf(\phi)-s\phi)\phi_{t}d\boldsymbol{x}\right).} \label{MCH3}
		\end{align}
	}
	\begin{remark}\leavevmode
		\begin{itemize}
			\item[(i)]
			 The right-hand side of \eqref{MCH3} vanishes at the continuous level, since $r(t)\equiv 1$. This allows us to introduce an arbitrarily small parameter $\alpha$, which can effectively slow down the variation of $r$ and keep it close to $1$, without altering the time scale of the dynamics of $\phi$. In the next section, we show that, {\color{black}under certain assumptions (including the boundedness assumption of $f^{\prime}(x)$)}, choosing a small $\alpha$ improves both the stability and accuracy of the numerical scheme.
			\item[(ii)] 
			{\color{black}When $\alpha$ tends to infinity, our RLM method recovers the original LM method.} However, in \cite{cheng2020}, one needs to solve a nonlinear algebraic equation for the Lagrange multiplier $q$.
			As we will see later, we can design a simple, explicit (for nonlinear term) linear scheme with constant coefficient for the reformulated system ($\ref{MCH1}$)-($\ref{MCH3}$).
		\end{itemize}
	\end{remark}
	\begin{lemma}\label{th1}
		The reformulated system $(\ref{MCH1})$-$(\ref{MCH3})$ obeys the following energy dissipation law:
		\begin{align}
			\frac{d}{dt}E_{RM}(\phi)=-\int_{\Omega}\mu\mathcal{G}\mu d\boldsymbol{x},\label{MED}
		\end{align}
		where 
		\begin{align}
			E_{RM}(\phi)=\frac{1}{2}\int_{\Omega}\phi\mathcal{L}\phi d\boldsymbol{x}+\int_{\Omega}F(\phi)d\boldsymbol{x}+{\frac{r-1}{\alpha}}.\label{ECM}
		\end{align}
	\end{lemma}
	\begin{proof} By taking the inner product of $(\ref{MCH1})$ with $\mu$ in the $L^{2}$ space and taking the $L^{2}$ inner product of ($\ref{MCH2}$) with $\phi_{t}$, we obtain
		\begin{align}
			\frac{1}{2}\frac{d}{dt}\int_{\Omega}\phi\mathcal{L}\phi d\boldsymbol{x}+\frac{s}{2}\frac{d}{dt}\int_{\Omega}|\phi|^{2}d\boldsymbol{x}+\int_{\Omega}(rf(\phi)-s\phi)\phi_{t}d\boldsymbol{x}=-\int_{\Omega}\mu\mathcal{G}\mu d\boldsymbol{x}.\label{CHME4}
		\end{align}
		By multiplying \eqref{MCH3} with $\frac{1}{\alpha}$ and combining it and $(\ref{CHME4})$, we get the energy dissipation law immediately
		\begin{align}
			\frac{1}{2}\frac{d}{dt}\int_{\Omega}\phi\mathcal{L}\phi d\boldsymbol{x}+\frac{d}{dt}\int_{\Omega}F(\phi)d\boldsymbol{x}+\frac{1}{\alpha}\frac{d(r-1)}{dt}=-\int_{\Omega}\mu\mathcal{G}\mu d\boldsymbol{x}\leq 0.
			\label{CHME8}
		\end{align} 
	\end{proof}
    {\color{black}
		\begin{remark}\leavevmode
			Notice that the extra term ${\frac{r-1}{\alpha}}$ in the energy functional $E_{RM}(\phi)$ is zero at continuous PDE level since $r\equiv 1$
		\end{remark}
	}
	\begin{remark}  It is worth noting that $E_{RM}(\phi)$ and the original energy $E(\phi)$ in equation \eqref{EO} are equivalent, { differing only by a constant $\frac{r-1}{\alpha}$ which is equal to $0$ on the continuous level. Moreoever, on the dicrete level, $r^n-1$ is small when $\alpha$ is chosen to be small. We will explain it in the next section.}
    \end{remark}
	\subsection{Numerical scheme}
	In this section, we construct a semi-implicit first-order backward Euler (BE) scheme, a second-order backward differentiation formula (BDF2) scheme and a second-order Crank-Nicolson (CN) scheme for the reformulated system \eqref{MCH1}--\eqref{MCH3}. Additionally, we provide a proof that all of these schemes are unconditionally energy stable.

	\subsection{The first-order RLM backward Euler (RLM-BE) scheme}
	We first present the following semi-implicit first-order RLM-BE scheme with the stabilization parameter $s>0$ and the relaxation parameter $\alpha>0$:
	\begin{align}
		\frac{\phi^{n+1}-\phi^{n}}{{\color{black}\tau}}&=-\mathcal{G}\mu^{n+1}, \label{NuSCH1}\\
		\mu^{n+1}&=\mathcal{L}\phi^{n+1}+s\phi^{n+1}+{(r^{n}f(\phi^{n})-s\phi^{n})}, \label{NuSCH2}\\
		\frac{r^{n+1}-r^{n}}{{\color{black}\tau}}=\alpha\bigg(-&\frac{E_{0}(\phi^{n+1})-E_{0}(\phi^{n})}{{\color{black}\tau}}+{\int_{\Omega}(r^{n}f(\phi^{n})-s\phi^{n})\frac{\phi^{n+1}-\phi^{n}}{{\color{black}\tau}}d\boldsymbol{x}}\bigg). \label{NuSCH3}
	\end{align}
	\begin{theorem}\label{th2} The first-order scheme ($\ref{NuSCH1}$)-($\ref{NuSCH3}$) for the equivalent system \eqref{MCH1}--\eqref{MCH3} is unconditionally energy stable in the sense that\\
		\begin{align}
			E^{n+1}_{RM}-E^{n}_{RM} =& -{\color{black}\tau}\int_{\Omega}\mu^{n+1}\mathcal{G}\mu^{n+1}d\boldsymbol{x}\nonumber\\
			&-\frac{1}{2}\int_{\Omega}(\phi^{n+1}-\phi^{n})\mathcal{L}(\phi^{n+1}-\phi^{n})d\boldsymbol{x}-\frac{s}{2}\int_{\Omega}|\phi^{n+1}-\phi^{n}|^{2}d\boldsymbol{x} \le 0, \label{NuED}
		\end{align}
		where the modified energy is defined by 
		\begin{align}
			E^{n}_{RM}=&\frac{1}{2}\int_{\Omega}\phi^{n}\mathcal{L}\phi^{n}d\boldsymbol{x}+\int_{\Omega}F(\phi^{n})d\boldsymbol{x}+{\frac{r^{n}-1}{\alpha}}.\label{NuED2}
		\end{align}
	\end{theorem}
	\begin{proof} We take the $L^{2}$ inner product of ($\ref{NuSCH1}$) with ${\color{black}\tau}\mu^{n+1}$ and take the $L^{2}$ inner product of ($\ref{NuSCH2}$) with $(\phi^{n+1}-\phi^{n})$. Then, we obtain
		\begin{align}
			&\frac{1}{2}\int_{\Omega}\big(\phi^{n+1}\mathcal{L}\phi^{n+1}-\phi^{n}\mathcal{L}\phi^{n}+(\phi^{n+1}-\phi^{n})\mathcal{L}(\phi^{n+1}-\phi^{n})\big)d\boldsymbol{x}+\frac{s}{2}\big(\int_{\Omega}|\phi^{n+1}|^{2}-|\phi^{n}|^{2}\nonumber\\
			&+|\phi^{n+1}-\phi^{n}|^{2}d\boldsymbol{x}\big)+\int_{\Omega}(r^{n}f(\phi^{n})-s\phi^{n})(\phi^{n+1}-\phi^{n})d\boldsymbol{x}=-{\color{black}\tau}\int_{\Omega}\mu^{n+1}\mathcal{G}\mu^{n+1}d\boldsymbol{x}. \label{NuCHME3}
		\end{align}
		By multiplying $(\ref{NuSCH3})$ with $\frac{{\color{black}\tau}}{\alpha}$ and combining it with (\ref{NuCHME3}), we can get
		\begin{align}
			&\frac{1}{2}\int_{\Omega}\big(\phi^{n+1}\mathcal{L}\phi^{n+1}-\phi^{n}\mathcal{L}\phi^{n}+(\phi^{n+1}-\phi^{n})\mathcal{L}(\phi^{n+1}-\phi^{n})\big)d\boldsymbol{x}+\frac{1}{\alpha}r^{n+1}-\frac{1}{\alpha}r^{n}+\nonumber\\
			&\int_{\Omega}F(\phi^{n+1})d\boldsymbol{x}-\int_{\Omega}F(\phi^{n})d\boldsymbol{x}+\frac{s}{2}\int_{\Omega}|\phi^{n+1}-\phi^{n}|^{2}d\boldsymbol{x}=-{\color{black}\tau}\int_{\Omega}\mu^{n+1}\mathcal{G}\mu^{n+1}d\boldsymbol{x}.\label{NuED}
		\end{align}
	\end{proof}
    {\color{black}
		\begin{remark}\leavevmode
			It is shown in Theorem \ref{th7} that, for the Allan-Cahn equation, $\frac{r^{n}-1}{\alpha}$ converges to $0$, as ${\color{black}\tau} \to 0$, which  means the difference between the discrete modified energy \eqref{NuED2} and the original energy is small.
		\end{remark}
	}
	
	\subsection{The second-order RLM Crank-Nicolson (RLM-CN) scheme}
	We consider the following second-order RLM-CN scheme with the stabilization parameter $s>0$ and the relaxation parameter $\alpha>0$:
	\begin{align}
		\frac{\phi^{n+1}-\phi^{n}}{{\color{black}\tau}}&=-\mathcal{G}\mu^{n+\frac{1}{2}}, \label{SeONuSCH1}\\
		\mu^{n+\frac{1}{2}}&=\mathcal{L} \phi^{n+\frac{1}{2}}+s\phi^{n+\frac{1}{2}}+{\left(\overline{r}^{n+\frac{1}{2}}f(\overline{\phi}^{n+\frac{1}{2}})-s\overline{\phi}^{n+\frac{1}{2}}\right)}, \label{SeONuSCH2}\\ 
		\frac{r^{n+1}-r^{n}}{{\color{black}\tau}}&=\alpha\bigg(-\frac{E_{0}(\phi^{n+1})-E_{0}(\phi^{n})}{{\color{black}\tau}}\nonumber\\
		&+{\int_{\Omega}\big(\overline{r}^{n+\frac{1}{2}}f(\overline{\phi}^{n+\frac{1}{2}})-s\overline{\phi}^{n+\frac{1}{2}}\big)\frac{\phi^{n+1}-\phi^{n}}{{\color{black}\tau}}d\boldsymbol{x}}\bigg),\label{SeONuSCH3}
	\end{align}
	where $\phi^{n+\frac{1}{2}}=\frac{1}{2}(\phi^{n+1}+\phi^{n})$, $\overline{\phi}^{n+\frac{1}{2}}=(\frac{3}{2}\phi^{n}-\frac{1}{2}\phi^{n-1})$, $\overline{r}^{n+\frac{1}{2}}=(\frac{3}{2}r^{n}-\frac{1}{2}r^{n-1})$.
	\begin{theorem}\label{th3} The second-order RLM-CN scheme ($\ref{SeONuSCH1}$)-($\ref{SeONuSCH3}$) for the equivalent system \eqref{MCH1}--\eqref{MCH3} is unconditionally energy stable in the sense that
		\begin{align}
			&E^{n+1}_{RM}-E^{n}_{RM} = -{\color{black}\tau}\int_{\Omega}\mu^{n+\frac{1}{2}}\mathcal{G}\mu^{n+\frac{1}{2}}d\boldsymbol{x} \le 0, \label{NuED}
		\end{align}
		where the modified energy is defined by 
		\begin{align}
			E^{n}_{RM}=&\frac{1}{2}\int_{\Omega}\phi^{n}\mathcal{L}\phi^{n}d\boldsymbol{x}+\int_{\Omega}F(\phi^{n})d\boldsymbol{x}+\frac{r^{n}-1}{\alpha}.\label{NuED_CN}
		\end{align}
	\end{theorem}
	\begin{proof} Similar to Theorem \ref{th2}.
	\end{proof}
	\subsection{The second-order RLM backward differentiation formula (RLM-BDF2) scheme}
	The second-order RLM-BDF2 scheme is given as 
	\begin{align}
		\frac{3 \phi^{n+1}-4 \phi^n+\phi^{n-1}}{2 {\color{black}\tau}}&=-\mathcal{G}\mu^{n+1}, \label{B2NuSCH1}\\
		\mu^{n+1}&=\mathcal{L} \phi^{n+1}+s\phi^{n+1}+(\overline{r}^{n+\frac{1}{2}}f(\overline{\phi}^{n+\frac{1}{2}})-s\overline{\phi}^{n+\frac{1}{2}}), \label{B2NuSCH2}\\
		\frac{3 r^{n+1}-4 r^n+r^{n-1}}{2 {\color{black}\tau}}&=\alpha\bigg(-\frac{3E_{0}(\phi^{n+1})-4E_{0}(\phi^{n})+E_{0}(\phi^{n-1})}{2{\color{black}\tau}} \nonumber\\
		&+\int_{\Omega}(\overline{r}^{n+\frac{1}{2}}f(\overline{\phi}^{n+\frac{1}{2}})-s\overline{\phi}^{n+\frac{1}{2}})\frac{3\phi^{n+1}-4\phi^{n}+\phi^{n-1}}{2{\color{black}\tau}}d\boldsymbol{x}\bigg),\label{B2NuSCH3}
	\end{align}
	where {$s>0$ is the stabilization parameter, $\alpha>0$ is the relaxation parameter}, $\overline{\phi}^{n+\frac{1}{2}}=(2\phi^{n}-\phi^{n-1})$ and $\overline{r}^{n+\frac{1}{2}}=(2 r^{n}-r^{n-1})$.
	\begin{theorem}\label{th4} The second-order RLM-BDF2 scheme ($\ref{B2NuSCH1}$)-($\ref{B2NuSCH3}$) for the equivalent system \eqref{MCH1}--\eqref{MCH3} is unconditionally energy stable in the sense that
		\begin{align}
			E^{n+1}_{RM}-E^{n}_{RM} =& -{\color{black}\tau}\int_{\Omega}\mu^{n+\frac{1}{2}}\mathcal{G}\mu^{n+\frac{1}{2}}d\boldsymbol{x} - \frac{1}{4}\int_{\Omega}(\phi^{n+1}-2\phi^{n}+\phi^{n-1})\mathcal{L}(\phi^{n+1}-2\phi^{n}\nonumber\\
			&+\phi^{n-1})d\boldsymbol{x} -\frac{s}{4}\int_{\Omega}|\phi^{n+1}-2\phi^{n}+\phi^{n-1}|^{2}\big)d\boldsymbol{x}\le 0, \label{NuED_BDF}
		\end{align}
		where the modified discrete version of the energy is defined by 
		\begin{align}
			&E^{n}_{RM}=\frac{1}{4}\int_{\Omega}\big(\phi^{n}\mathcal{L}\phi^{n}+(2\phi^{n}-\phi^{n-1})\mathcal{L}(2\phi^{n}-\phi^{n-1})\big)d\boldsymbol{x}+\frac{s}{4}\int_{\Omega}\big(|\phi^{n}|^{2}+ \nonumber\\
			&|2\phi^{n}-\phi^{n-1}|^{2}\big)d\boldsymbol{x}+\frac{1}{\alpha}(\frac{1}{2}(3r^{n}-r^{n-1})-1)+\frac{1}{2}\big(3E_{0}(\phi^{n})-E_{0}(\phi^{n-1})\big).\label{NuED2_BDF}
		\end{align}
	\end{theorem}
	\begin{proof} We take the $L^{2}$ inner product of ($\ref{B2NuSCH1}$) with ${\color{black}\tau}\mu^{n+1}$ and take the $L^{2}$ inner product of ($\ref{B2NuSCH2}$) with $\frac{1}{2}(3 \phi^{n+1}-4 \phi^n+\phi^{n-1})$. Then, we can get
		\begin{align}
			&\frac{1}{4}\int_{\Omega}\big(\phi^{n+1}\mathcal{L}\phi^{n+1}+(2\phi^{n+1}-\phi^{n})\mathcal{L}(2\phi^{n+1}-\phi^{n})-\phi^{n}\mathcal{L}\phi^{n}-\nonumber\\
			&(2\phi^{n}-\phi^{n-1})\mathcal{L}(2\phi^{n}-\phi^{n-1})+(\phi^{n+1}-2\phi^{n}+\phi^{n-1})\mathcal{L}(\phi^{n+1}-2\phi^{n}+\phi^{n-1})\big)d\boldsymbol{x}+ \nonumber\\
			&\frac{s}{4}\int_{\Omega}\big(|\phi^{n+1}|^{2}+|2\phi^{n+1}-\phi^{n}|^{2}-|\phi^{n}|^{2}-|2\phi^{n}-\phi^{n-1}|^{2}+|\phi^{n+1}-2\phi^{n}+\phi^{n-1}|^{2}\big)d\boldsymbol{x} \nonumber\\
			&+\frac{1}{2}\int_{\Omega}(\overline{r}^{n+\frac{1}{2}}f(\overline{\phi}^{n+\frac{1}{2}})-s\overline{\phi}^{n+\frac{1}{2}})(3 \phi^{n+1}-4 \phi^n+\phi^{n-1})d\boldsymbol{x}=-{\color{black}\tau}\int_{\Omega}\mu^{n+1}\mathcal{G}\mu^{n+1}d\boldsymbol{x}. \label{B2NuCHME3}
		\end{align}
		By multiplying $(\ref{B2NuSCH3})$ with $\frac{{\color{black}\tau}}{\alpha}$ and combining it with (\ref{B2NuCHME3})
		\begin{align}
			&\frac{1}{4}\int_{\Omega}\big(\phi^{n+1}\mathcal{L}\phi^{n+1}+(2\phi^{n+1}-\phi^{n})\mathcal{L}(2\phi^{n+1}-\phi^{n})-\phi^{n}\mathcal{L}\phi^{n}-\nonumber\\
			&(2\phi^{n}-\phi^{n-1})\mathcal{L}(2\phi^{n}-\phi^{n-1})+(\phi^{n+1}-2\phi^{n}+\phi^{n-1})\mathcal{L}(\phi^{n+1}-2\phi^{n}+\phi^{n-1})\big)d\boldsymbol{x}+ \nonumber\\
			&\frac{s}{4}\int_{\Omega}\big(|\phi^{n+1}|^{2}+|2\phi^{n+1}-\phi^{n}|^{2}-|\phi^{n}|^{2}-|2\phi^{n}-\phi^{n-1}|^{2}+|\phi^{n+1}-2\phi^{n}+\phi^{n-1}|^{2}\big)d\boldsymbol{x} \nonumber\\
			&+\frac{1}{2\alpha}(3r^{n+1}-r^n)-\frac{1}{2\alpha}(3r^{n}-r^{n-1})+\frac{1}{2}\big(3E_{0}(\phi^{n+1})-E_{0}(\phi^{n})\big)-\nonumber\\
			&\frac{1}{2}\big(3E_{0}(\phi^{n})-E_{0}(\phi^{n-1})\big) =-{\color{black}\tau}\int_{\Omega}\mu^{n+1}\mathcal{G}\mu^{n+1}d\boldsymbol{x}. \label{B2NuCHME5}
		\end{align}
	\end{proof}
	\section{The boundedness of the Lagrange multiplier $r^n$}
	We now demonstrate the boundedness of the Lagrange multiplier $r$ using the Allen-Cahn equation and the Cahn-Hilliard equation as examples. In this proof, we assume that the second derivative of $F(\phi)$ is bounded. {\color{black}Consequently, we adjust $F(\phi)$ to exhibit a quadratic growth rate for $|\phi|>1$, which is a commonly adopted practice \cite{shen2010,yu2017numerical}:
\begin{align}
{F}(\phi)=\frac{1}{4}
\begin{cases}
4(\phi+1)^2, & \text{if } \phi<-1, \\
(\phi^2-1)^2, & \text{if } -1 \leq \phi \leq 1, \\
4(\phi-1)^2, & \text{if } \phi>1.
\end{cases}
\label{trun_F}
\end{align}}
Accordingly, $f(\phi)$ is defined as $f(\phi)=F'(\phi)$. Throughout the paper, $\|\cdot\|_{V}$ denotes the norm in the space $V$, while $\|\cdot\|$ denotes the $L^{2}$ norm.

		\subsection{Allen-Cahn equation}
		{\color{black} According to the numerical schemes  \eqref{NuSCH1}--\eqref{NuSCH3} and \eqref{B2NuSCH1}--\eqref{B2NuSCH3}}, we can rewrite the first-order RLM-BE scheme and the second-order RLM-BDF2 scheme for Allen-Cahn equation as:
		\begin{align}
			& \frac{\phi^{n+1}-\phi^{n}}{{\color{black}\tau}}=\Delta\phi^{n+1}-r^{n}f(\phi^{n})-s(\phi^{n+1}-\phi^{n}), \label{AC_1st_1}\\
			& \frac{r^{n+1}-r^{n}}{{\color{black}\tau}}=\alpha\bigg(-\frac{E_{0}(\phi^{n+1})-E_{0}(\phi^{n})}{{\color{black}\tau}}+{\int_{\Omega}(r^{n}f(\phi^{n})-s\phi^{n})\frac{\phi^{n+1}-\phi^{n}}{{\color{black}\tau}}d\boldsymbol{x}}\bigg), \label{AC_1st_2}
		\end{align}
		and
		\begin{align}
			\frac{3\phi^{n+1}-4\phi^{n}+\phi^{n-1}}{2{\color{black}\tau}}=&\Delta\phi^{n+1}-\color{black}\overline{r}^{n+\frac{1}{2}}f(\overline{\phi}^{n+\frac{1}{2}})-\color{black}s(\phi^{n+1}-\overline{\phi}^{n+\frac{1}{2}}), \label{AC_nu1}\\
			\frac{3r^{n+1}-4r^{n}+r^{n-1}}{2{\color{black}\tau}}=&\alpha\big(-\frac{3E^{n+1}_{0}-4E^{n}_{0}+E^{n-1}_{0}}{2{\color{black}\tau}}+\color{black}\int_{\Omega}(\overline{r}^{n+\frac{1}{2}}f(\overline{\phi}^{n+\frac{1}{2}})-s\overline{\phi}^{n+\frac{1}{2}})\nonumber\\
			&\frac{3\phi^{n+1}-4\phi^{n}+\phi^{n-1}}{2{\color{black}\tau}}d\boldsymbol{x}\big). \label{AC_nu2}
		\end{align}
		\begin{theorem}\label{th5}
			Consider the numerical schemes \eqref{AC_1st_1}--\eqref{AC_1st_2} and \eqref{AC_nu1}--\eqref{AC_nu2}. Assume: (i) $\phi^0\in H^{3}(\Omega)$ with $\Omega\subseteq\mathbb {R}^d$, $d=1,2,3$;
			(ii) $|f^{\prime}(x)|<M_{0}$, $\forall x\in \mathbb {R}$. 
			For any given $M_{1}>1$,  there exists $ \gamma_{0}>0$ depending on $T$, $s$, $\phi^{0}$, $\phi^{1}$, $M_{0}$, $M_{1}$, and $\Omega$, such that $\forall \alpha<\gamma_0 {\color{black}\tau}$, $r^{m}\in [0,M_{1}]$ and $\|\phi^{m}\|_{H^{2}}\le M$, where the constant $M$ depends on $T$, $s$, $\phi^{0}$, $\phi^{1}$, $M_{0}$, $M_{1}$, and $\Omega$.
		\end{theorem}
		\begin{proof} Because of the similarity between the proofs of first-order and second-order numerical schemes, for brevity, we only present the proof of the second-order scheme in \ref{AP2}.
		\end{proof}
		
		\subsection{Cahn-Hilliard equation}
		Similarly, we can rewrite the first-order RLM-BE scheme and the second-order RLM-BDF2 scheme for Cahn-Hilliard equation as:
		\begin{align}
			& \frac{\phi^{n+1}-\phi^{n}}{{\color{black}\tau}}=-\Delta^{2}\phi^{n+1}+r^{n}\Delta f(\phi^{n})+s\Delta(\phi^{n+1}-\phi^{n}), \label{CH_1st_1}\\
			& \frac{r^{n+1}-r^{n}}{{\color{black}\tau}}=\alpha\bigg(-\frac{E_{0}(\phi^{n+1})-E_{0}(\phi^{n})}{{\color{black}\tau}}+{\int_{\Omega}(r^{n}f(\phi^{n})-s\phi^{n})\frac{\phi^{n+1}-\phi^{n}}{{\color{black}\tau}}d\boldsymbol{x}}\bigg), \label{CH_1st_2}
		\end{align}
		and
		\begin{align}
			\frac{3\phi^{n+1}-4\phi^{n}+\phi^{n-1}}{2{\color{black}\tau}}=&-\Delta^{2}\phi^{n+1}+\color{black}\overline{r}^{n+\frac{1}{2}}\Delta f(\overline{\phi}^{n+\frac{1}{2}})+s\Delta(\phi^{n+1}-\overline{\phi}^{n+\frac{1}{2}}), \label{CH_nu1}\\
			\frac{3r^{n+1}-4r^{n}+r^{n-1}}{2{\color{black}\tau}}=&\alpha\big(-\frac{3E^{n+1}_{0}-4E^{n}_{0}+E^{n-1}_{0}}{2{\color{black}\tau}}+\color{black}\int_{\Omega}(\overline{r}^{n+\frac{1}{2}}f(\overline{\phi}^{n+\frac{1}{2}})-\nonumber\\
			&\color{black}s\overline{\phi}^{n+\frac{1}{2}})\frac{3\phi^{n+1}-4\phi^{n}+\phi^{n-1}}{2{\color{black}\tau}}d\boldsymbol{x}\big) . \label{CH_nu2}
		\end{align}
		\begin{theorem}\label{th6}
			Consider the numerical schemes \eqref{CH_1st_1}--\eqref{CH_1st_2} and \eqref{CH_nu1}--\eqref{CH_nu2}.
			Assume:
			(i) $\phi^0\in H^{4}(\Omega)$ with $\Omega\subseteq\mathbb {R}^d$, $d=1,2,3$;
			(ii) $|f^{\prime}(x)|<M_{0}$, $\forall x\in \mathbb {R}$. 
			For any given $M_{1}>1$,  there exists $ \gamma_{0}>0$ depending on $T$, $s$, $\phi^{0}$, $\phi^{1}$, $M_{0}$, $M_{1}$, and $\Omega$, such that $\forall \alpha<\gamma_0 {\color{black}\tau}$, $r^{m}\in [0,M_{1}]$ and $\|\phi^{m}\|_{H^{2}}\le M$, where the constant $M$ depends on $T$, $s$, $\phi^{0}$, $\phi^{1}$, $M_{0}$, $M_{1}$, and $\Omega$.
		\end{theorem}
		\begin{proof} For the same reason as Theorem \ref{th5}, we only present the proof of the second-order scheme in \ref{AP3}.\end{proof}
		\begin{remark} 
        In Example 3 of Section \ref{Nuexp}, our numerical experiments show that the truncation of double-well potential function \eqref{trun_F} guarantees the boundedness of the Lagrange multiplier $r^n$, enhancing the numerical stability of the algorithm.
        \end{remark}
		
	
	\section{Further discussions on the first-order RLM scheme}
	It is well-known that the Allen-Cahn equation not only {\color{black}satisfies} the energy dissipation law but also the maximum principle. In particular, the first-order scheme usually preserves this property. This inspires us to establish the maximum principle of the first-order RLM scheme \eqref{AC_1st_1}--\eqref{AC_1st_2}. Moreover, we study the behavior of the added term $\frac{r^n-1}{\alpha}$ appearing in modified energy \eqref{NuED2} when ${\color{black}\tau}$ and $\alpha$ tends to $0$. In addition, we give an error estimate where the dependency on $\alpha$ is given explicitly. 
	
	Let $\phi(t)=\phi(x,t)$ be the solution of the original system  \eqref{GF}, and $\phi^n=\phi^n(x)$ and $r^{n}$ be the numerical solution of the reformulated system \eqref{MCH1}--\eqref{MCH3}. Denote $z^{n+1}=|r^{n+1}-1|$ and  $e^{n+1}=\phi^{n+1}-\phi(t^{n+1})$.

	

	\begin{theorem}[Maximum principle]\label{AP1th2} Assume that $\|\phi^{0}\|_{\infty}\le 1$, $r^{0}\in [0,M_1]$ with $M_1>1$ some constant and $s\ge M_1\|f^{\prime}(x)\|_{{C}[-1,1]}$. Then, there exists $\gamma_{0}$ depending on $s$ and $M_1$, such that $\forall\alpha\leq {\color{black}\tau}\gamma_{0}$,  the first-order scheme \eqref{AC_1st_1}--\eqref{AC_1st_2} satisfies $r^{m}\in [0,M_1]$ and $ \|\phi^{m}\|_{\infty}\le 1$.
	\end{theorem}
	\begin{proof} Assume $r^{n}\in [0,M_1]$ and $\|\phi^{n}\|_{\infty}\le 1$ for all $n\leq m-1$. 
		From \eqref{AC_1st_1}, we can get 
		\begin{align}
			\phi^{n+1}=\left[\left(\frac{1}{{\color{black}\tau}}+s\right)\mathcal{I}-\Delta\right]^{-1}\left(\frac{1}{{\color{black}\tau}}\phi^{n}+\left(r^{n}\left(-f(\phi^{n})\right)+s\phi^{n}\right)\right),
		\end{align}
		where $\mathcal I$ is the identity operator.
		By Lemma \ref{lm2} in \ref{AP1}, we have 
		\begin{align}
		\left\|\left[\left(\frac{1}{{\color{black}\tau}}+s\right)\mathcal{I}-\Delta\right]^{-1}\right\|_{\infty}\le \left(\frac{1}{{\color{black}\tau}}+s\right)^{-1}.
		\end{align}
		Since $s\ge M_1\|f^{\prime}(x)\|_{{C}[-1,1]}$, $r^{n}\in [0,M_1]$, and $\|\phi^{n}\|_{\infty}\le1$, it follows from Lemma~\ref{lm1} in Appendix~\ref{AP1} that
		\begin{align}
			\left\|\left(\frac{1}{{\color{black}\tau}}\phi^{n}+\left(r^{n}\left(-f(\phi^{n})\right)+s\phi^{n}\right)\right)\right\|_{\infty}&\le\frac{1}{{\color{black}\tau}}\|\phi^{n}\|_{\infty}+\|\left(r^{n}\left(-f(\phi^{n})\right)+s\phi^{n}\right)\|_{\infty}\nonumber\\
			&\le\left(\frac{1}{{\color{black}\tau}}+s\right).
		\end{align}
		Therefore, we obtain $\|\phi^{n+1}\|_{\infty}\le1$, implying $\|\phi^{m}\|_{\infty}\le1$.
		
		Next, we prove $r^{m}\in [0,M_1]$. Let $z^{n}=|r^{n}-1|$, from \eqref{AC_1st_2}, we can get 
		\begin{align}
			&z^{n+1}-z^{n}=\left(|r^{n+1}-1|-|r^{n}-1|\right)\le |(r^{n+1}-1)-(r^{n}-1)|\nonumber\\
			&\le\alpha\left|-(E_{0}(\phi^{n+1})-E_{0}(\phi^{n}))+\int_{\Omega}\left(f(\phi^{n})-s\phi^{n}\right)(\phi^{n+1}-\phi^{n})d\boldsymbol{x}\right|\nonumber\\
			&+\alpha z^{n}\int_{\Omega}\left|f(\phi^{n})(\phi^{n+1}-\phi^{n})\right|d\boldsymbol{x} \le\alpha C+\alpha z^{n}C, \label{zn1}
		\end{align}
		where we have used the fact that $E_0$ and $f$ are polynomials.
		Taking the summation of \eqref{zn1} for $n$ from $0$ to $m-1$,
		\begin{align}
			z^{m}\le mC\alpha + \sum_{n=0}^{m-1}\alpha C z^{n}.
		\end{align}
		By setting $\alpha={\color{black}\tau}\gamma$, then we can get
		$
		z^{m}\le TC\gamma + {\color{black}\tau}\sum_{n=0}^{m-1}\gamma C z^{n}.
		$
		By applying Lemma \ref{lm3} in \ref{AP1}, we can get
		\begin{align}
			\color{black} z^{m}\le e^{TC\gamma}TC\gamma,
		\end{align}
		We can get $ \lim_{\gamma \to 0}z^{m}=0$, so there exists a $\gamma_{0}$, such that $z^{m}<\min(1,M_1-1)$. Therefore, we can get $r^{m}\in [0,M_1]$.
	\end{proof}
	
	\begin{theorem}\label{th7}
		{Assume that $\phi(t_{0})\in {H}^{3}(\Omega)$, $\phi_{t}\in{L}^{\infty}(0,T;{L}^{2}(\Omega))$, $\phi_{tt}\in{L}^{2}(0,T;{L}^{2}(\Omega))$, $\|\phi(t_{0})\|_{\infty}\le 1$, $s\ge M\|f
			^{\prime}(x)\|_{{C}[-1,1]}$, with some constant $M>1$. {\color{black}For the first-order scheme \eqref{AC_1st_1}--\eqref{AC_1st_2}}, we then have for all $n\le T/{\color{black}\tau}$, } 
		\begin{itemize}
			\item[(i)] $ \lim_{\alpha \to 0}|r^{n}-1|=0$;
			\item[(ii)]   $\lim_{{\color{black}\tau}\to 0}\frac{1}{\alpha}|r^{n}-1|=0$;\
			\item[(iii)] 
			There exists $\gamma_{1} > 0$, depending on $T$, $s$, $M$, $\phi(t_{0})$, and $\Omega$, such that $\forall \alpha \le \gamma_{1}\tau$,
			\begin{align}
				\frac{1}{2}||\nabla e^{n}||^{2}+\frac{s}{2}||e^{n}||^{2}\le C \exp\left((1-C{\color{black}\tau})^{-1}Ct^{n}\right)&\big({\color{black}\tau}^{2}\int_{0}^{t_{n}}\left(||\phi_{tt}(s)||^{2}+||\phi_{t}(s)||^{2}\right)ds\nonumber\\
                &+{\tau^2}\gamma_{1}\big), \nonumber
			\end{align}
			where the constant $C$ depends on $T$, $s$, $M$, $\phi(t_{0})$, $\Omega$, $||\phi||_{{L}^{\infty}(0,T;W^{1,\infty}(\Omega))}$, and $||\phi_{t}||_{{L}^{\infty}(0,T;{L}^{2}(\Omega))}$.
		\end{itemize}
	\end{theorem}
	\begin{proof} See \ref{AP4} for detailed proof of (i) and (ii). For the proof of (iii), one can follow the framework in \cite{shen2018convergence} by deriving first the $H^{2}$ bound of the solution and then the error bound. For brevity, we omit the details here.
	\end{proof}

	
	The results show that under certain conditions, the solution of equations \eqref{AC_1st_1}--\eqref{AC_1st_2} converges to the solution of the original model, and $\frac{r^{n}-1}{\alpha}$ converges to $0$, as ${\color{black}\tau}$ approach to zero,  which also means the difference between the discrete modified energy \eqref{NuED2} and the original energy is small.
	
	In \cite{zhang2025afully} and in the later numerical examples (Section \ref{AC and CH results}), we will observe that the numerical solutions achieve higher accuracy as $\alpha$ decreases.

	\section{Numerical results}\label{Nuexp}
	In this section, we now implement the proposed numerical methods and apply them to several classical phase-field models, including the Allen-Cahn (AC) equation, the Cahn-Hilliard (CH) equation and the molecular beam epitaxy (MBE) model.
	
	For this study, we consider phase-field models with periodic boundary conditions. The finite element method (FEM) is employed for spatial discretization. The spatial domain $\Omega=[0, L_x]\times[0, L_y]$ is divided into a uniform grid with a mesh size of $h_x=L_x/N_x$ and $h_y=L_y/N_y$, where $N_x$ and $N_y$ are positive even integers. Unless otherwise stated, we set $N_x=N_y=128$ for all cases. Since both the BDF2 and CN schemes offer second-order accuracy, our focus in this paper is on comparing the RLM-CN scheme with the LM-CN scheme, the original SAV-CN scheme and the RSAV-CN scheme.
	{The reference  energy values in the following numerical examples are obtained using the LM-CN schemes with a very small time step size ${\color{black}\tau}=10^{-5}$ to ensure the uniqueness of the solution to nonlinear algebraic equations.} The RLM modified energy is equal to $E_{RM}$ in equation \eqref{ECM} after subtracting a constant $\frac{r-1}{\alpha}$.
	\par
	\subsection{Allen-Cahn and Cahn-Hilliard equations}\label{AC and CH results}
    {\color{black}Consider the free energy
    $
    E(\phi)=\int_\Omega \frac{\epsilon^2}{2}|\nabla\phi|^2+\frac14(\phi^2-1)^2\,dx.
    $
    When the mobility operator is $\mathcal G=\lambda>0$, the general phase field model in \eqref{GF} blackuces to the Allen--Cahn (AC) equation.}
	If  the  mobility operator is $\mathcal{G}=-\lambda\Delta$, the phase field model in \eqref{GF} blackuces to the Cahn-Hilliard (CH) equation.
	
\begin{figure}
		\centering
		\subfigure[Convergence tests]{
			\label{Fig.sub.1}
			\includegraphics[width=55mm]{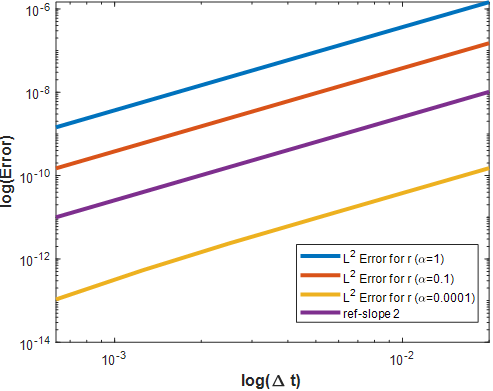}}
		\hspace{2.5cm}
		\subfigure[Convergence tests]{
			\label{Fig.sub.2}
			\includegraphics[width=55mm]{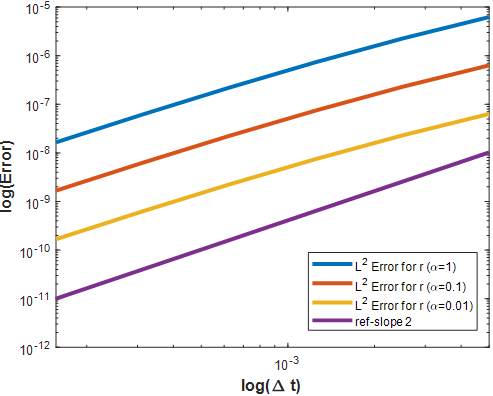}}
		\caption{(a) Convergence tests of the Lagrange multiplier $r$ at t=0.5 for different $\alpha$ in the Allen-Cahn equation. (b) Convergence tests of the Lagrange multiplier $r$ at t=0.125 for different $\alpha$ in the Cahn-Hilliard equation.}
		\label{fig:Fig.1}
	\end{figure}

	\textbf{Example 1.} We first verify that the RLM-CN scheme is second-order accurate in time. For the Allen-Cahn equation, consider the domain $\Omega=[0,2\pi]^{2}$, and we choose the smooth initial condition
	{
		\begin{align}
			\phi_0(x, y, 0)=\sum_{i=1}^2\tanh \left(\frac{\sqrt{\left(x-x_i\right)^2+\left(y-y_i\right)^2}-R_i}{\sqrt{2} \epsilon}\right)+1, \label{AC_ex1_init1}
		\end{align}
		where $R_{1}=1.4$, $R_{2}=0.5$, ${(x_{1},y_{1})}=(\pi-0.8,\pi)$, and ${(x_{2},y_{2})}=(\pi+1.7,\pi)$.
		To ensure that the initial value of the phase field variable satisfies periodic boundary conditions, we make the following modification to the initial condition \eqref{AC_ex1_init1}:
		\begin{align}
			\phi_0 = 
			\begin{cases}
				\sum_{i=1}^2\tanh \left(\frac{\sqrt{\left(x-x_i\right)^2+\left(y-y_i\right)^2}-R_i}{\sqrt{2} \epsilon}\right)+1, & \text{if}\quad (x,y)\in [\delta,2\pi-\delta]\times[\delta,2\pi-\delta], \\
				-1, & \text{otherwise},
			\end{cases}
			\label{AC_ex1_init}
		\end{align}
		where $\delta=0.1$.} To solve the Allen-Cahn equation, we use uniform meshes $N_x=N_y=256$ and set the model parameters $\epsilon=0.08$ and $\lambda=1$, and the numerical parameter $s=4$. We choose three different $\alpha=1, \ 0.1, \ 0.0001$. Since the analytical solutions are unavailable, we use the numerical solutions obtained with a smaller time step size and smaller $\alpha$, namely ${\color{black}\tau}=3.125\times 10^{-4}$ and $\alpha=1\times 10^{-5}$, as a reference for  the exact solution. The numerical errors of $\phi$ in $L^{2}$ norm with different time steps and different $\alpha$ at $t=0.5$ are summarized in Table~\ref{err_time_AC}, and the numerical errors of $r$ are summarized in the Fig.\ref{Fig.sub.1}. We observe second-order convergence for both the numerical solutions of $\phi$ and $r$, and the accuracy of $\phi$ and $r$ increases as the value of $\alpha$ decreases, which also verifies the conclusion of Theorem \ref{th7}. { In terms of the computational efficiency, we show in Table \ref{comp_eff_AC} a comparison of compute time using SAV method, LM method and RLM method for solving this example up to $t=0.5$, both with second-order temporal accuracy. Notably, the RLM method is slightly more efficient, requiring roughly half the computational time of the SAV and LM methods.}\par
	For the Cahn-Hilliard equation, we use the same initial condition  $\eqref{AC_ex1_init}$. The model parameters are set to $\epsilon=0.08$ and $\lambda=0.125$. We solve the CH equation using uniform meshes $N_x=N_y=256$, with numerical parameters $s=4$ and $\alpha=1,0.1,0.01$. Since the analytical solutions are unknown, we consider the numerical solutions obtained with a smaller time step size and smaller $\alpha$, namely ${{\color{black}\tau}}=7.8125\times 10^{-5}$ and $\alpha=1\times 10^{-5}$, as an approximation of the exact solution. The numerical errors of $\phi$ in $L^{2}$ norm with different time steps at $t=0.125$ are summarized in Table~\ref{err_time_CH}, and the numerical errors of $r$ are summarized in the Fig.\ref{Fig.sub.2}. We observe second-order convergence for both the numerical solutions of $\phi$ and $r$, and the accuracy of $\phi$ and $r$ increases as the relaxation parameter $\alpha$ decreases.\par
\begin{figure}
		\centering
		\subfigure[Stability tests.]{
			\label{Fig2.sub.1}
			\includegraphics[width=52mm]{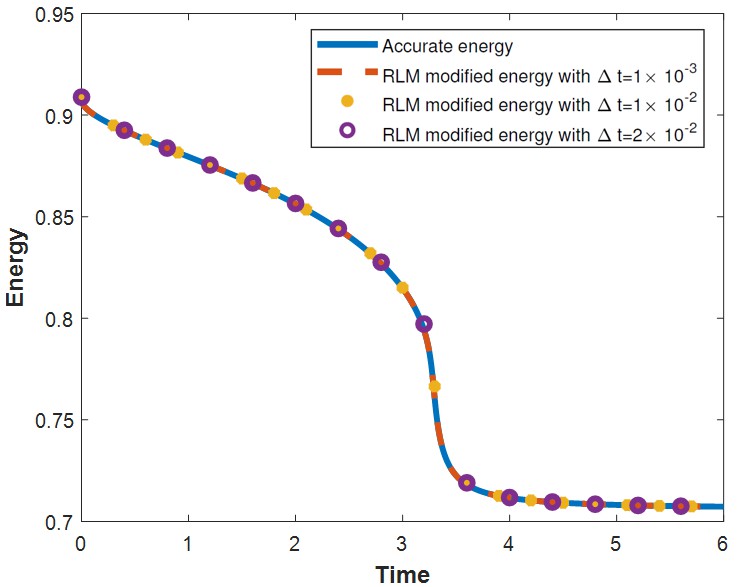}}
		\hspace{2.5cm}
		\subfigure[The evolution of energy functional for different $\alpha$]{
			\label{Fig2.sub.2}
			\includegraphics[width=51mm]{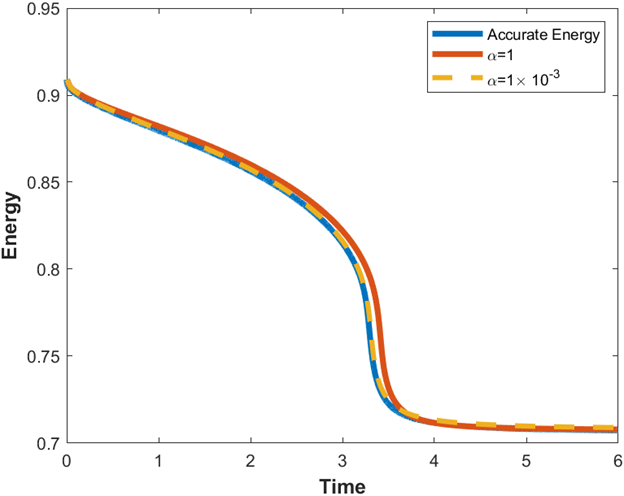}}
		\caption{(a) The time evolution of the energy functional for the different time step sizes ${\color{black}\tau}=2\times 10^{-2}$, $1\times 10^{-2}$, $1\times 10^{-3}$ in Example 2. (b) The time evolution of the energy functional for the different $\alpha$ in Example 2.}
		\label{fig:Fig.2}
	\end{figure}

\begin{figure}
		\centering
		\subfigure[The evolution of $r^{n}$]{
			\label{Fig3.sub.1}
			\includegraphics[width=45mm]{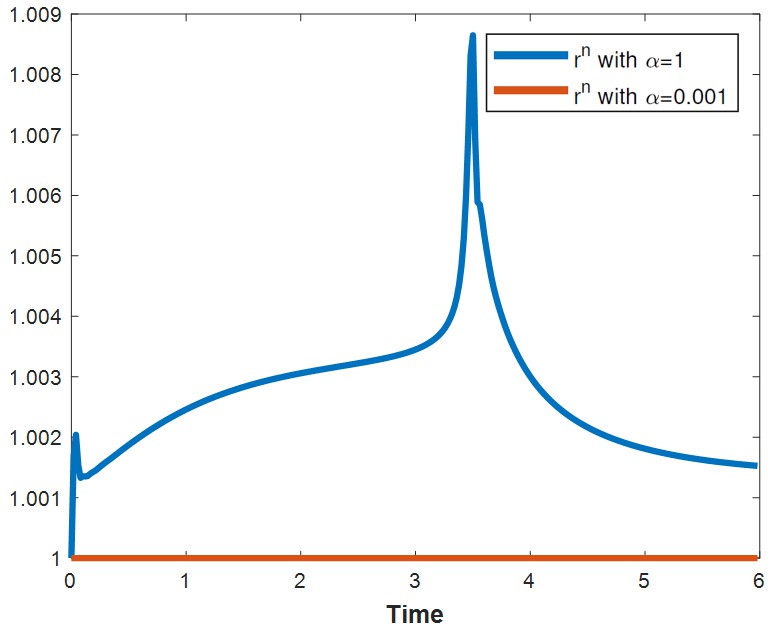}}
		\hspace{2.5cm}
		\subfigure[The evolution of $\phi$]{
			\label{Fig3.sub.2}
			\includegraphics[width=40mm]{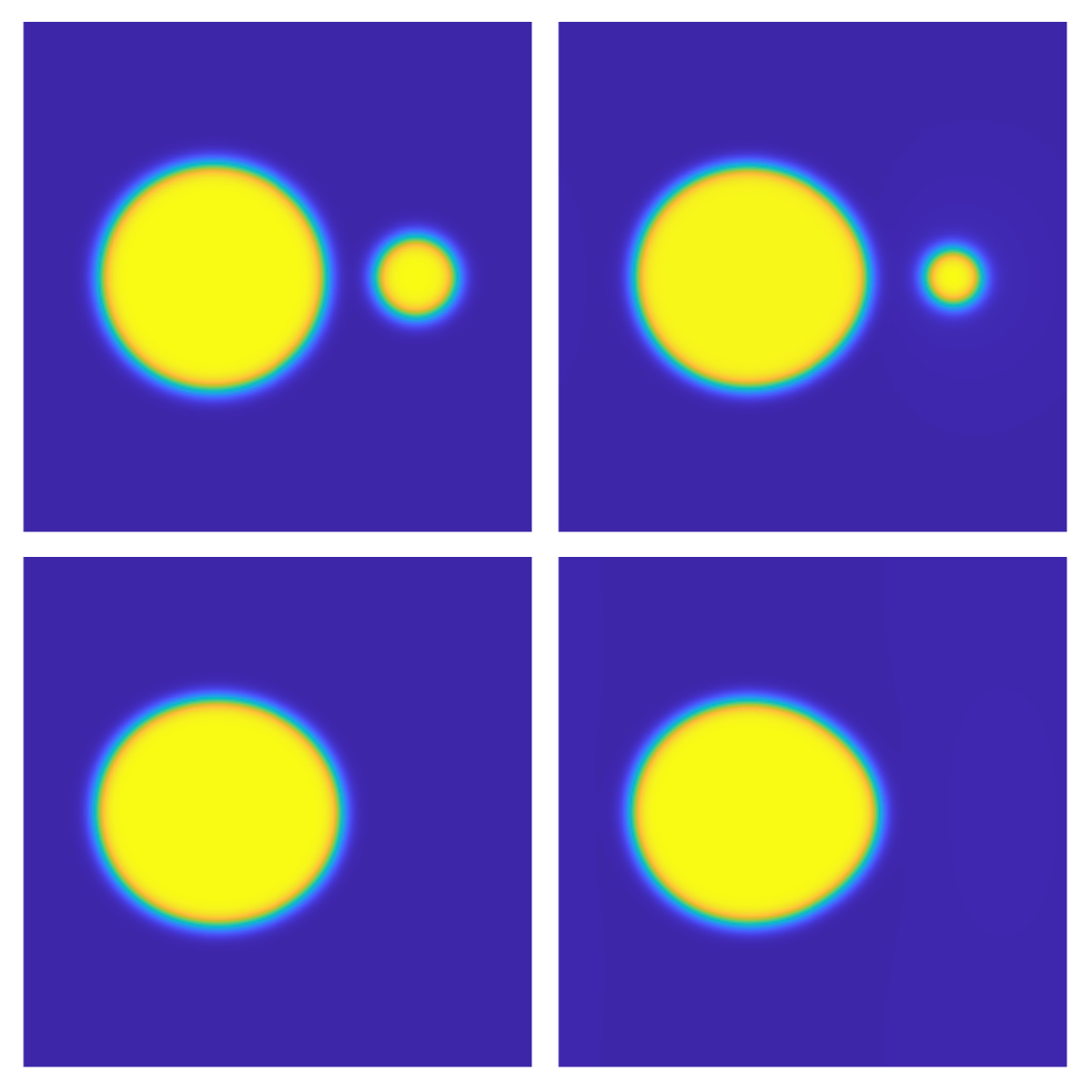}}
		\caption{(a) The evolution of $r^{n}$ for different $\alpha$ in Example 2. (b) Snapshots of the phase variable $\phi$ are taken at t=0, 2, 4, 6 in clockwise sense in Example 2.}
		\label{fig:Fig.3}
	\end{figure}

	\begin{table}\footnotesize
		\begin{center}
			\begin{tabular}{|c|c|c|c|c|c|c|}
				\hline ${\color{black}\tau}$ & $\|e(\phi)\|_{L^2}$ ($\alpha=1$) & Order & $\|e(\phi)\|_{L^2}$ ($\alpha=0.1$) & Order & $\|e(\phi)\|_{L^2}$ ($\alpha=10^{-4}$) & Order \\
				\hline $6.25\times 10^{-4}$ & $1.405\times 10^{-9}$ & & $9.844\times 10^{-10}$ & & $9.823\times 10^{-10}$ & \\
				\hline $1.25\times 10^{-3}$ & $6.314\times 10^{-9}$ & 2.17 &  $4.914\times 10^{-9}$ & 2.32 & $4.912\times 10^{-9}$ & 2.32 \\
				\hline $2.5\times 10^{-3}$ & $2.597\times 10^{-8}$ & 2.04 & $2.065\times 10^{-8}$ & 2.07 & $2.063\times 10^{-8}$ & 2.07  \\
				\hline $5\times 10^{-3}$ & $1.043\times 10^{-7}$ & 2.01 & $8.352\times 10^{-8}$ & 2.02 & $8.349\times 10^{-8}$ & 2.02 \\
				\hline $1\times 10^{-2}$ & $4.157\times 10^{-7}$ & 1.99 & $3.350\times 10^{-7}$ & 2.00 & $3.349\times 10^{-7}$ & 2.00  \\
				\hline $2\times 10^{-2}$ & $1.646\times 10^{-6}$ & 1.99 & $1.340\times 10^{-6}$ & 2.00 & $1.340\times 10^{-6}$ & 2.00 \\
				\hline
			\end{tabular}
			\caption{\label{err_time_AC}{Convergence tests where the $L^{2}$ numerical errors of the phase-field variable $\phi$ at $t=0.5$ for different $\alpha=1,0.1,0.0001$ in the Allen-Cahn equation.} }
		\end{center}
	\end{table}
	{
		\begin{table}\footnotesize
			{
				\color{black}
				\begin{center}
					\begin{tabular}{|c|lc|lc|lc|}
						\hline & SAV & & RLM & & LM &\\
						\hline ${\color{black}\tau}$ & $\|e(\phi)\|_{L^2}$ & CPU time(s) & $\|e(\phi)\|_{L^2}$ & CPU time(s) & $\|e(\phi)\|_{L^2}$ & CPU time(s) \\
						\hline $2.5\times 10^{-3}$ & $4.525\times 10^{-8}$ & 31.37 & $4.511\times 10^{-8}$ & 16.50 & $4.513\times 10^{-8}$ & 32.79  \\
						\hline $5\times 10^{-3}$ & $1.829\times 10^{-7}$ & 15.41 & $1.826\times 10^{-7}$ & 7.81 & $1.826\times 10^{-7}$ & 16.56 \\
						\hline $1\times 10^{-2}$ & $7.342\times 10^{-7}$ & 7.79 & $7.324\times 10^{-7}$ & 4.03 & $7.325\times 10^{-7}$ & 8.29  \\
						\hline $2\times 10^{-2}$ & $2.945\times 10^{-6}$ & 3.96 & $2.931\times 10^{-6}$ & 1.98 & $2.933\times 10^{-6}$ & 4.21 \\
						\hline
					\end{tabular}
					\caption{\label{comp_eff_AC}{Comparison of computational cost of the SAV method, LM method, and RLM method for the Allen-Cahn equation.} }
				\end{center}
			}
		\end{table}
	}
	\begin{table}\footnotesize
		\begin{center}
			\begin{tabular}{|c|c|c|c|c|c|c|}
				\hline ${\color{black}\tau}$ & $\|e(\phi)\|_{L^2}$ ($\alpha=1$) & Order & $\|e(\phi)\|_{L^2}$ ($\alpha=0.1$) & Order & $\|e(\phi)\|_{L^2}$ ($\alpha=0.01$) & Order \\
				\hline $1.5625\times 10^{-4}$ & $1.09\times 10^{-8}$ & & $5.43\times 10^{-9}$ & & $5.03\times 10^{-9}$ & \\
				\hline $3.125\times 10^{-4}$ & $4.60\times 10^{-8}$ & 2.07 & $2.66\times 10^{-8}$ & 2.29 & $2.51\times 10^{-8}$ & 2.32 \\
				\hline $6.25\times 10^{-4}$ & $1.78\times 10^{-7}$ & 1.95 & $1.11\times 10^{-7}$ & 2.06 & $1.06\times 10^{-7}$ & 2.08  \\
				\hline $1.25\times 10^{-3}$ & $6.65\times 10^{-7}$ & 1.90 & $4.46\times 10^{-7}$ & 2.01 & $4.29\times 10^{-7}$ & 2.02 \\
				\hline $2.5\times 10^{-3}$ & $2.43\times 10^{-6}$ & 1.87 & $1.78\times 10^{-6}$ & 2.00 & $1.73\times 10^{-6}$ & 2.01  \\
				\hline $5\times 10^{-3}$ & $8.88\times 10^{-6}$ & 1.87 & $7.18\times 10^{-6}$ & 2.01 & $7.03\times 10^{-6}$ & 2.02 \\
				\hline
			\end{tabular}
			\caption{\label{err_time_CH}{Convergence tests where the $L^{2}$ numerical errors of the phase-field variable $\phi$ at $t=0.125$ for different $\alpha=1,0.1,0.01$ in the Cahn-Hilliard equation.} }
		\end{center}
	\end{table}
	\textbf{Example 2.} In this example, we perform the energy stability tests and investigate the impact of different $\alpha$ on the numerical results. Taking the Cahn-Hilliard equation as an example, we consider the domain $\Omega=[0,2\pi]^{2}$, and use the same initial condition as the \textbf{Example 1}. 
	We set the model parameters $\epsilon=0.08$ and $\lambda=1$ and numerical parameters $s=4$. In Fig. \ref{Fig2.sub.1}, we set $\alpha=10^{-3}$ and plot the time evolution of the energy functional for the different time step sizes ${\color{black}\tau}=2\times 10^{-2}$, $1\times 10^{-2}$, $1\times 10^{-3}$. We observe that all energy curves decay monotonically for all time steps, confirming the unconditional energy stability of the RLM scheme. In Fig. \ref{Fig2.sub.2}, we set ${\color{black}\tau}=1\times 10^{-2}$ and plot the time evolution of the energy functional {\color{black} $E_{RM}-\frac{r-1}{\alpha}$} for different values of $\alpha$: $\alpha=1$ and $1\times 10^{-3}$. We observe that the value of $\alpha$ does not affect energy dissipation, and as $\alpha$ decreases, the system's energy approaches the accurate energy. In Fig. \ref{Fig3.sub.1}, we show the evolution of $r^{n}$ for different $\alpha$, $\alpha=1$, $1\times 10^{-3}$. In Figure \ref{Fig3.sub.2}, we illustrate the evolution of $\phi$ over time. As time evolves, the fluid within smaller bubbles dissolves into the surrounding fluid and diffuses towards larger bubbles, causing them to grow. This observation agrees with the phenomenon of Ostwald ripening in an open system. The accurate simulation demonstrates the effectiveness of our RLM approach.
	\\

\begin{figure}
		\centering
		\subfigure[The original energy evolutions when ${\color{black}\tau}=0.05$.]{
			\label{Fig5.sub.1}
			\includegraphics[width=40mm]{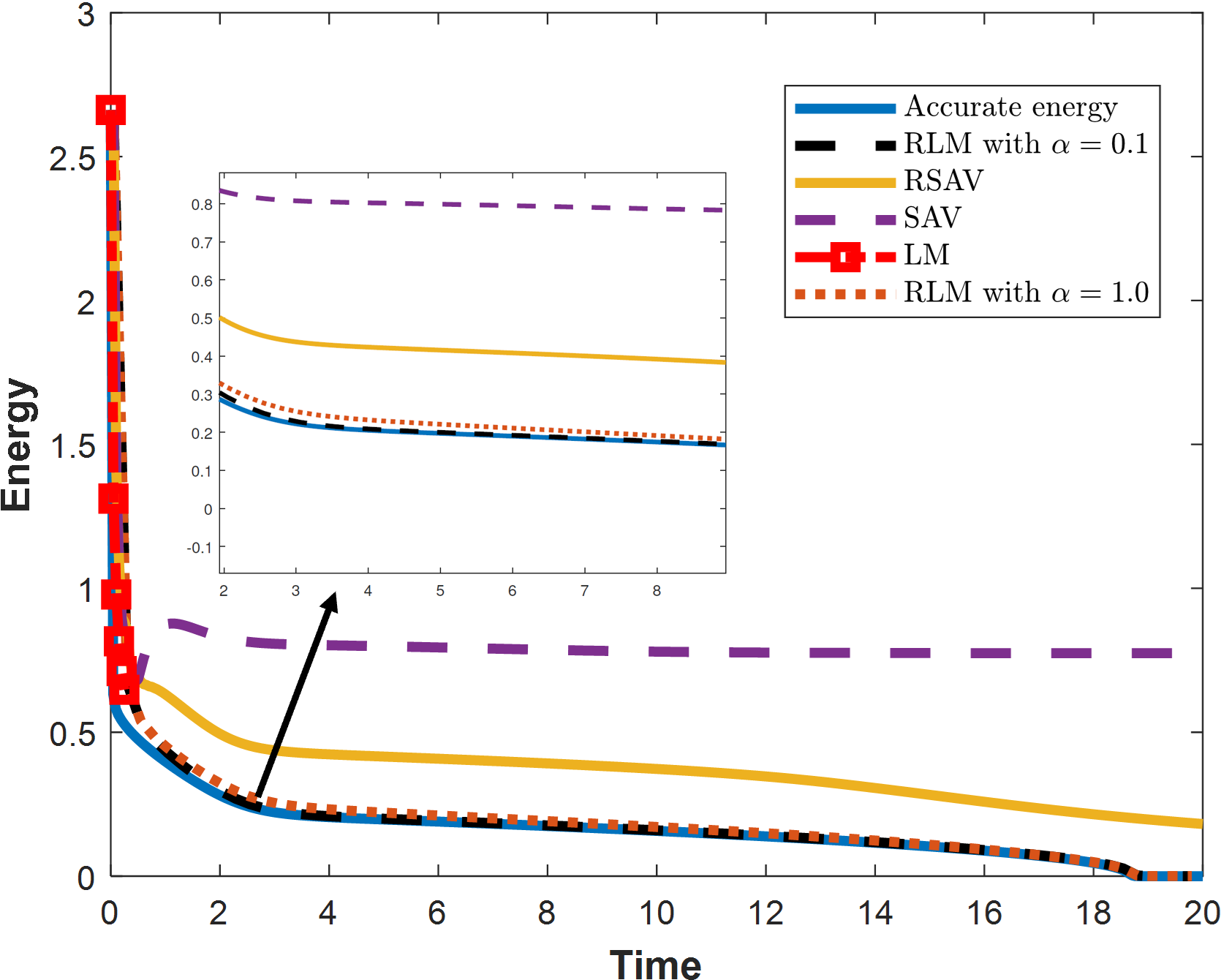}}
		\hspace{0.5cm}    
		\subfigure[The original energy evolutions when ${\color{black}\tau}=0.025$ for the SAV, RSAV, RLM methods and ${\color{black}\tau}=0.0125$ for the LM method.]{
			\label{Fig5.sub.3}
			\includegraphics[width=40mm]{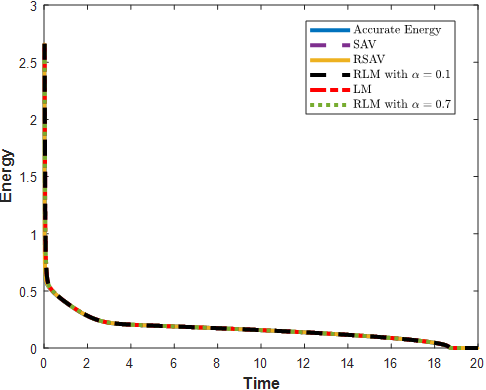}}    
		\hspace{0.5cm}
		\subfigure[The comparison of $\frac{\eta}{\sqrt{E_0(\overline{\phi}) + C_0}}$ in SAV, $q^{n+\frac{1}{2}}$ in LM, and $\overline{r}^{n+\frac{1}{2}}$ in RLM when ${\color{black}\tau}=0.05$.]{
			\label{Fig5.sub.2}
			\includegraphics[width=39mm]{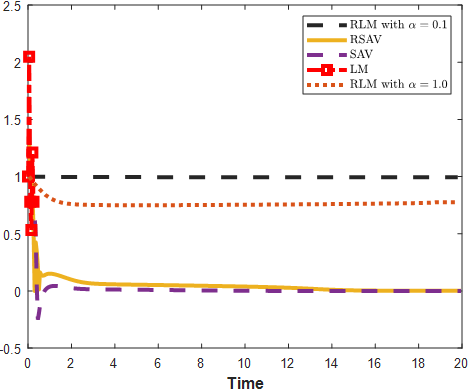}}
		\caption{(a) and (b) The comparisons of the LM original energy, the RLM modified energy, the RSAV modified energy, the SAV modified energy for the Allen-Cahn equation when using various time step sizes in Example 3. (c) The comparisons of $\frac{\eta}{\sqrt{E_0(\overline{\phi}) + C_0}}$ in the SAV-CN and RSAV-CN methods, $q^{n+\frac{1}{2}}$ in \eqref{eq:LM_nu_c2}, and $\overline{r}^{n+\frac{1}{2}}$ in \eqref{SeONuSCH2} for the Allen-Cahn equation in Example 3.}
		\label{fig:Fig.5}
	\end{figure}

\begin{figure}
		\centering 
		\centering 
		\includegraphics[width=30mm]{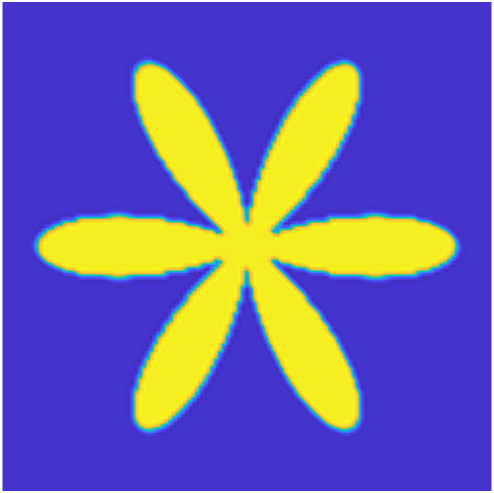} 
		\includegraphics[width=30mm]{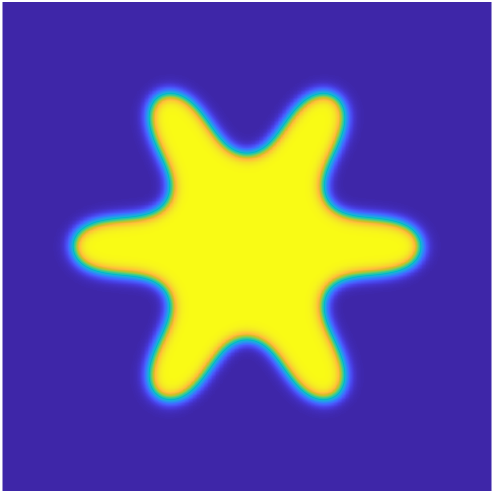}
		\includegraphics[width=30mm]{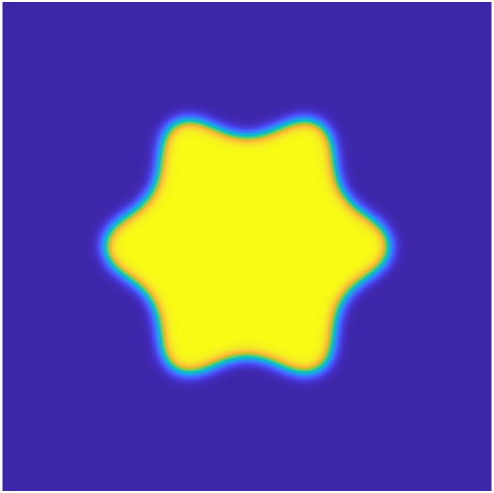}
		\includegraphics[width=30mm]{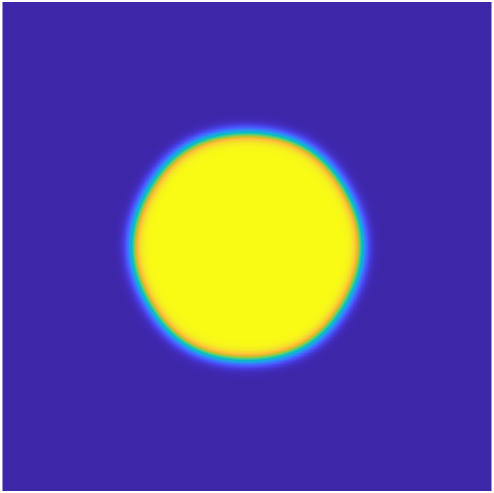}
		\caption{Snapshots of the phase variable $\phi$  at t=0, 2.5, 5, 10 for Example 3 (Allen-Cahn equation).}
		\label{fig:Fig.12}
	\end{figure}

	\textbf{Example 3.} In this example, we compare the performance of the RLM-CN scheme with the original SAV-CN scheme, {\color{black}the} RSAV-CN scheme and the LM-CN scheme. For the Allen-Cahn equation, we consider the domain $\Omega=[0,L_{x}]\times[0,L_{y}]$ and choose the initial condition as follows:
	\begin{align}
		& \phi(x, y)=\tanh \frac{1.5+1.2 \cos (6 \theta)-2 \pi r}{\sqrt{2} \epsilon} ,\label{AC_ex2_init1}\\
		& \theta=\arctan \frac{y-0.5 L_y}{x-0.5 L_x}, \quad r=\sqrt{\left(x-\frac{L_x}{2}\right)^2+\left(y-\frac{L_y}{2}\right)^2}. \label{AC_ex2_init2}
	\end{align}
	{\color{black}We consider the parameters $L_x=L_y=1$, $\epsilon = 0.015$, $\lambda = 10$, $s=0$, $C_{0}=0$ in \eqref{eq:EM}, ${\color{black}\tau}=0.05$. 
		Figure~\ref{fig:Fig.5} presents the evolution of the original energy computed by different methods, as well as a comparison of the quantities $\frac{\eta^{n+\frac{1}{2}}}{\sqrt{E_0(\overline{\phi}^{n+\frac{1}{2}})+C_0}}$ in \eqref{eq:SAV_d2} obtained from the SAV-CN and RSAV-CN methods, $q^{n+\frac{1}{2}}$ in the LM method  \eqref{eq:LM_nu_c2}, and $\overline{r}^{n+\frac{1}{2}}$ in the RLM method \eqref{SeONuSCH2} to assess the consistency between the original and modified models. Figure~\ref{Fig5.sub.1} demonstrates that the RLM-CN method achieves better accuracy and stability than both the RSAV-CN and SAV-CN methods when using the large time step size ${\color{black}\tau}=0.05$. Moreover, using a smaller parameter $\alpha = 0.1$ leads to better accuracy compablack to the case of $\alpha = 1.0$. As shown in Figure~\ref{Fig5.sub.2}, the multiplier $r$ in the RLM method remains close to 1 when $\alpha = 0.1$, indicating better consistency between the original and discrete models than the case of $\alpha = 1.0$. In the LM method, solving the nonlinear equation for the Lagrange multiplier $q$ fails to converge at the 6th time step, leading to a breakdown of the solution. For the SAV method, the original energy increases at certain time steps, and $\frac{\eta^{n+\frac{1}{2}}}{\sqrt{E_0(\overline{\phi}^{n+\frac{1}{2}})+C_0}}$  deviates significantly from 1. Although the RSAV method maintains stability of the original energy, the same deviation in $\frac{\eta^{n+\frac{1}{2}}}{\sqrt{E_0(\overline{\phi}^{n+\frac{1}{2}})+C_0}}$  from 1 causes noticeable differences between the computed and accurate energy. When we use smaller time steps, from Figures \ref{Fig5.sub.3}, we observe that the energy evolution curves of the LM method, RSAV method, SAV method, and our RLM method closely match the accurate energy.} {\color{black} In Figure \ref{fig:Fig.12}, we illustrate the evolution of $\phi^{n}$ over time. As time evolves, the interface first shrinks into a small circular shape and eventually disappears.}\par
\begin{figure}
		\centering
		\subfigure[The evolution of original energy when ${\color{black}\tau}=0.05$.]{
			\label{new_sub_1}
			\includegraphics[width=40mm]{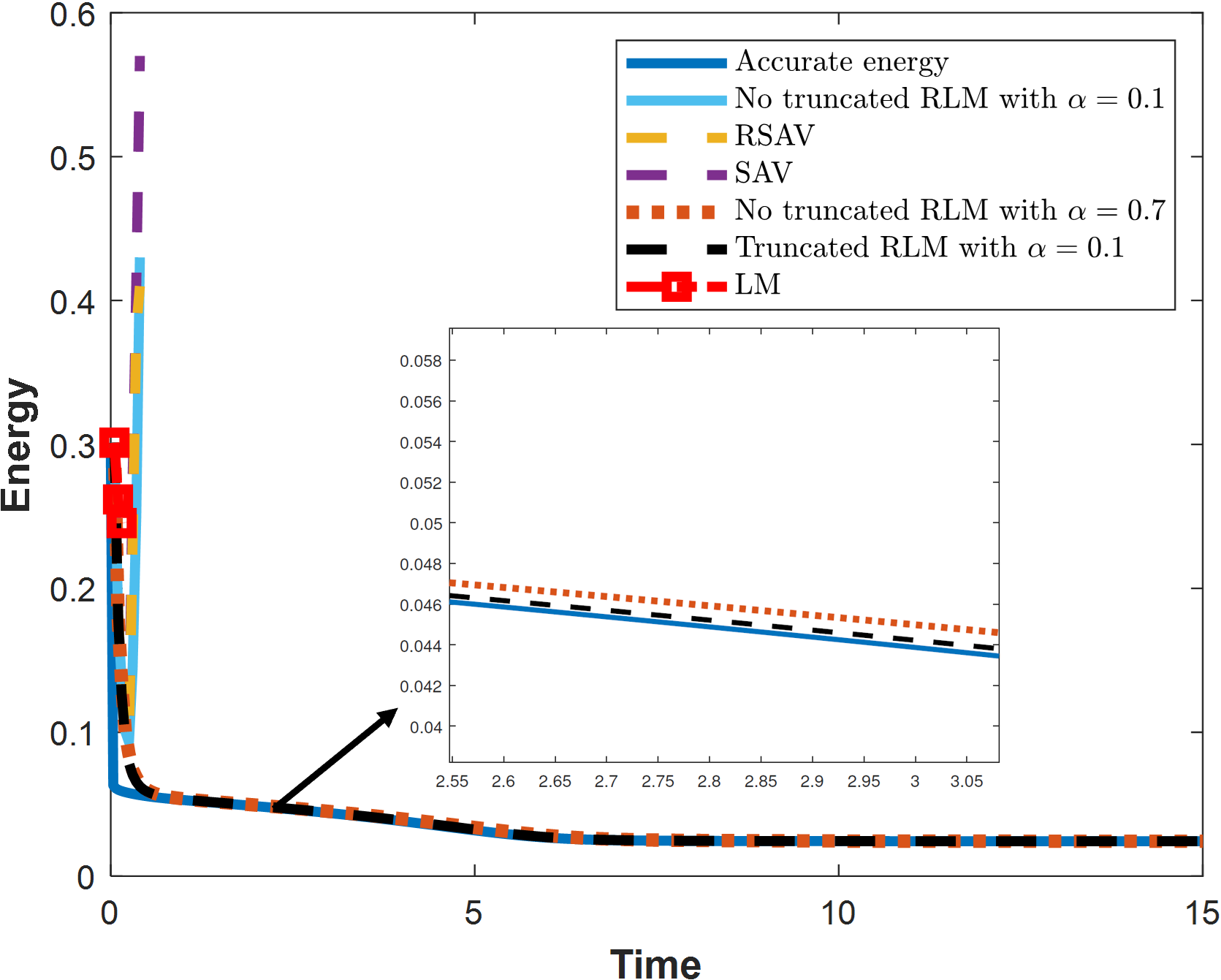}}
		\hspace{0.5cm}
		\subfigure[The evolution of original energy when ${\color{black}\tau}=0.01$ for the SAV, RSAV, RLM methods and ${\color{black}\tau}=0.001$ for the LM method.]{
			\label{new_sub_3}
			\includegraphics[width=40mm]{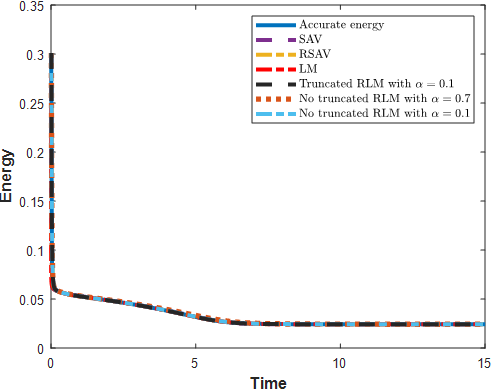}}  
		\hspace{0.5cm}
		\subfigure[The comparison of $\frac{\eta}{\sqrt{E_0(\overline{\phi}) + C_0}}$, $q^{n+\frac{1}{2}}$, and $\overline{r}^{n+\frac{1}{2}}$ when ${\color{black}\tau}=0.05$.]{
			\label{new_sub_2}
			\includegraphics[width=39mm]{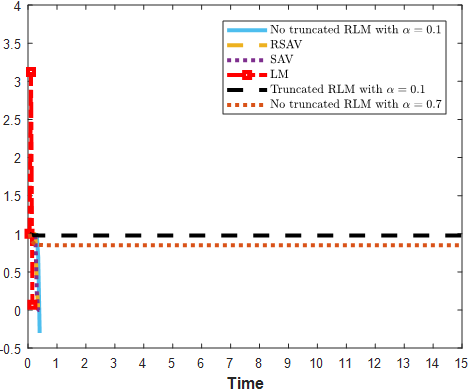}}
		\caption{(a) and (b) The comparisons of the LM original energy, the RLM original energy, the RSAV original energy, the SAV original energy for the Cahn-Hilliard equation when using some different time step sizes in Example 3. (c) The comparisons of $\frac{\eta}{\sqrt{E_0(\overline{\phi}) + C_0}}$ in the SAV-CN and RSAV-CN methods, $q^{n+\frac{1}{2}}$ in \eqref{eq:LM_nu_c2}, and $\overline{r}^{n+\frac{1}{2}}$ in \eqref{SeONuSCH2} for the Cahn-Hilliard equation in Example 3. No truncated RLM refers to the RLM method with the standard double well potential function, and truncated RLM refers to the RLM method with the truncated double well potential function \eqref{trun_F}.}
		\label{new_csav}
	\end{figure}

	{\color{black}For the Cahn-Hilliard equation, we choose the same initial condition and domain as in the above example. The model parameters used are $\lambda=0.01$, $\epsilon=0.015$, and we set the numerical parameters $s=3$, $C_{0}=10$ in \eqref{eq:EM} and ${\color{black}\tau}=0.05$. To validate the effectiveness of the truncated double-well potential in ensuring the boundedness of the Lagrange multiplier $r^n$ and to align with the statement of Theorem~\ref{th6}, we use the second-order RLM-BDF2 scheme \eqref{CH_nu1}--\eqref{CH_nu2}. Similar to the Allen--Cahn case, Figure \ref{new_csav} presents the evolution of the original energy obtained from different methods and a comparison of the quantities $\frac{\eta^{n+\frac{1}{2}}}{\sqrt{E_0(\overline{\phi}^{n+\frac{1}{2}})+C_0}}$ computed by the SAV-BDF2 and RSAV-BDF2 methods, $q^{n+\frac{1}{2}}$ in the LM method \eqref{eq:LM_nu_c2}, and $\overline{r}^{n+\frac{1}{2}}$ in the RLM method \eqref{B2NuSCH2}, to evaluate the consistency between the original and discrete models. Figure \ref{new_sub_1} demonstrates that the RLM-BDF2 method achieves better accuracy and stability compablack to the RSAV-BDF2 and SAV-BDF2 methods when ${\color{black}\tau}=0.05$. Moreover it is observed that when using a smaller parameter $\alpha = 0.1$ with the standard double-well potential, the original energy becomes unstable. From Figure \ref{new_sub_2}, this instability might be caused by the significant deviation of Lagrange multiplier from 1. Employing the truncated double-well potential \eqref{trun_F} can effectively eliminate this issue, keeping the Lagrange multiplier $r$ close to 1, which is also in agreement with the theoretical result stated in Theorem~\ref{th6}. Another way to avoid the energy instability in the RLM-BDF2 scheme is to choose a larger value of $\alpha$, such as $\alpha = 0.7$. In this case, the Lagrange multiplier $r$ slightly deviates from 1, which leads to a mild loss in accuracy, but on the other hand enhances the stability. Of course, this problem can be blackuced by taking a smaller time step (e.g., ${\color{black}\tau} = 0.01$), as observed in Figure \ref{new_sub_3}.}
	\subsection{Molecular beam epitaxy model with slope selection}
\begin{figure}
		\centering
		\subfigure[The energy evolutions]{
			\label{Fig7.sub.1}
			\includegraphics[width=40mm]{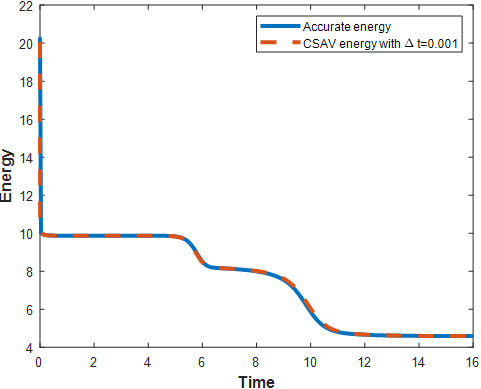}}
		\hspace{2.5cm}
		\subfigure[The evolution of $r^{n}$.]{
			\label{Fig7.sub.2}
			\includegraphics[width=43mm]{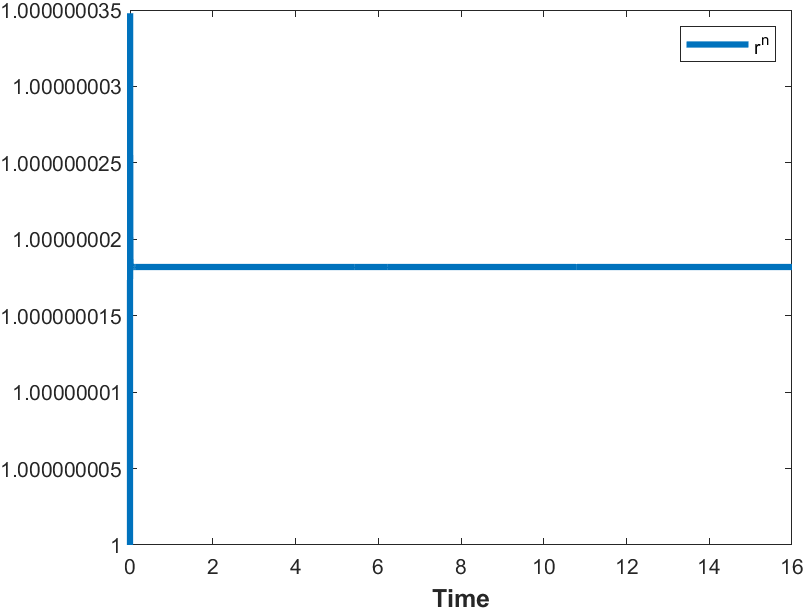}}
		\caption{(a) The evolution of energy computed by using the time step ${\color{black}\tau}=0.001$ for Example 4. (b) The evolution of $r^{n}$ for Example 4. }
		\label{fig:Fig.7}
	\end{figure}

\begin{figure}
		\centering 
		\centering 
		\includegraphics[width=30mm]{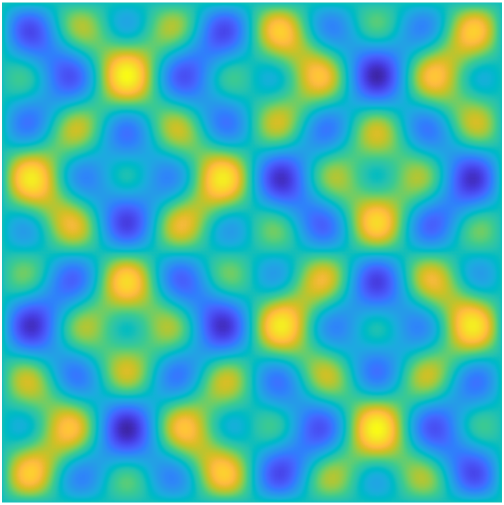} 
		\includegraphics[width=30mm]{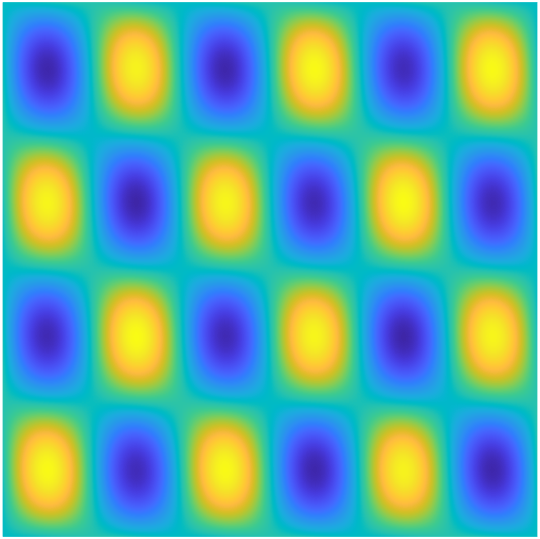}
		\includegraphics[width=30mm]{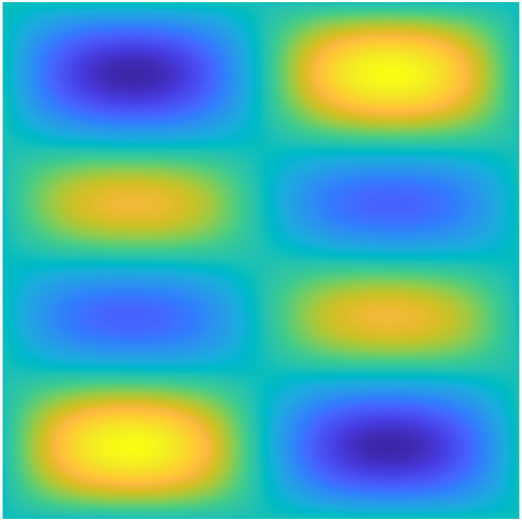}
		\includegraphics[width=30mm]{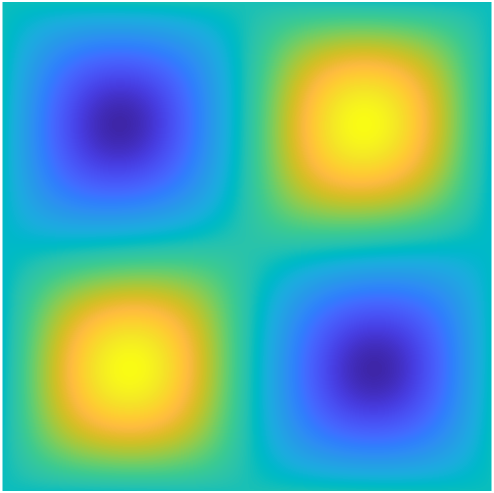}
		\caption{Snapshots of the phase variable $\phi$ at t=0, 1, 6, 16 for Example 4.}
		\label{fig:Fig.17}
	\end{figure}

	We investigate the molecular beam epitaxy (MBE) model with slope selection. $\phi$ represents the MBE thickness, and the free energy is defined as $E(\phi) =$\\ $ \int_{\Omega} \left(\frac{\epsilon^{2}}{2}(\Delta\phi)^{2} + \frac{1}{4}(\phi^{2}-1)^{2}\right)d\boldsymbol{x}$. The mobility operator is $\mathcal{G} = 1$.\\
	\textbf{Example 4.} We consider the domain $\Omega=[0,2\pi]^{2}$, with $\epsilon^{2}=0.1$.  we choose the smooth initial condition
	\begin{align}
		\phi(x,y,t=0)=0.1(\sin(3x)\sin(2y)+\sin(5x)\sin(5y)),\label{MBE_ex1_init}
	\end{align}
	and solve the MBE model with slope selection using $\alpha=10^{-5}$. In Fig.\ref{Fig7.sub.1}, we present the energy curve obtained with a time step of ${\color{black}\tau}=0.001$, which matches closely with the accurate energy. The evolution of $r^n$ is depicted in Fig.\ref{Fig7.sub.2}. It can be seen that $r^n$ remains  sufficiently close to 1, indicating a high numerical consistency between the modified and original systems. The solution contours at t = 0, 1, 6, and 16 are plotted in Fig.\ref{fig:Fig.17}, which agree with the results reported in \cite{xu2006}.
	\section{Conclusions}
	
	In this paper, we have proposed a relaxed Lagrange multiplier (RLM) approach that effectively addresses the inconsistency issue between the original model and the modified model, a problem that can arise in the original SAV method. We introduced a relaxation parameter $\alpha$ in the ODE for the Lagrange multiplier to slow its evolution, ensuring consistency between the original system and the reformulated system after numerical discretization. Furthermore, we rigorously prove that the numerical Lagrange multiplier $r^n$ has lower and upper bounds when $\alpha$ is small, and will tend to $1$ as $\alpha$ tends to $0$. Extensive numerical simulations have been conducted to demonstrate the accuracy and stability of the proposed method. 

    \section{Acknowledgments}
    C. Quan is supported by the National Natural Science Foundation of China (Grant No. 12271241), Guangdong Basic and Applied Basic Research Foundation (Grant No. 2023B1515020030),  Shenzhen Science and Technology Innovation Program (Grant No. JCYJ20230807092402004), and Hetao Shenzhen-Hong Kong Science and Technology  Innovation Cooperation Zone Project (No. HZQSWS-KCCYB-2024016). X.-P. Wang acknowledges support from the National Natural Science Foundation of China (NSFC) (No. 12271461), the key project of NSFC (No. 12131010), Shenzhen Science and Technology Innovation Program (Grant: C10120230046), the Hetao Shenzhen-Hong Kong Science and Technology  Innovation Cooperation Zone Project (No.HZQSWS-KCCYB-2024016) and the University
Development Fund from The Chinese University of Hong Kong, Shenzhen (UDF01002028).
	
	\appendix
	\section{Some fundamental lemmas and theorems}\label{AP1}
	{
		\noindent
		\begin{lemma}\label{lm1}
			Assume there exists a constant $\beta>0$ such that $Q(\beta)\le 0\le Q(-\beta)$. If $s\ge M\|Q^{\prime}(x)\|_{{C}[-\beta,\beta]}$, where $M$ is a positive constant, we have $|kQ(\xi)+s \xi| \leq s \beta$ for any $\xi \in[-\beta, \beta]$ and $k \in[0,M]$.
		\end{lemma}
		\begin{proof} Let $h(\xi)=kQ(\xi)+s\xi$. It follows that $h^{\prime}(\xi)=s+kQ^{\prime}(\xi)$, and 
			\begin{align}
				0\le h^{\prime}(\xi) \le 2s, \forall \xi\in[-\beta, \beta].
			\end{align}
			From the assumption, for any  $\xi\in[-\beta, \beta]$, we know that
			\begin{align}
				-s\beta\le -s\beta+kQ(-\beta)=h(-\beta)\le h(\xi)\le h(\beta)=s\beta+kQ(\beta)\le s\beta.
			\end{align}
		\end{proof}
		\begin{lemma}\label{lm2} For any $a>0$, we have $\|\left(a\mathcal{I}-\Delta\right)^{-1}\|_{\infty}\le a^{-1}$, where $\mathcal{I}$ represents the identity operator.
		\end{lemma}
		\begin{proof} Consider the equation \( a u - \Delta u = f \), where \( f \in L^\infty(\mathbb{R}^n) \). Suppose that \( u \) attains its positive maximum at a point \( x_0 \). At this point, we have
			\[
			\Delta u(x_0) \leq 0 \implies a u(x_0) \leq f(x_0) \leq \| f \|_\infty.
			\]
			Thus, \( u(x_0) \leq \| f \|_\infty / a \). Similarly, if \( u \) attains its negative minimum at \( x_1 \), then
			\[
			\Delta u(x_1) \geq 0 \implies a u(x_1) \geq f(x_1) \geq -\| f \|_\infty \implies |u(x_1)| \leq \| f \|_\infty / a.
			\]
			Combining these results, we obtain \( \| u \|_\infty \leq a^{-1} \| f \|_\infty \). Therefore, the operator norm satisfies
			\[
			\| (a \mathcal I - \Delta)^{-1} \|_{\infty} \leq a^{-1}.
			\]
		\end{proof}
		\begin{lemma}[Discrete Gr\"onwall Lemma]\label{lm3} Let $\left(u_n\right)$ and $\left(w_n\right)$ be nonnegative sequences satisfying
			$$
			u_n \leq \alpha+\sum_{k=0}^{n-1} u_k w_k, \quad \forall n.
			$$
			Then for all $n$ it holds
			$$
			u_n \leq \alpha e^{\sum_{k=0}^{n-1} w_k}.
			$$
		\end{lemma}
		\begin{lemma}\label{lm_shen}[Lemma 2.3 in \cite{shen2018convergence}]\label{lm4} (i) If $\|\phi\|_{H^1} \leq M_{0}$ and  $|g^{\prime}(x)|<C(|x|^{p}+1)$ $\big(p > 0\text { arbitrary if }$ $ n = 1, 2; 0 < p < 4 \text { if } n = 3\big)$, then  for any $\phi \in H^3(\Omega)$, there exist $0 \leq \sigma<1$ and a constant $C(M_{0})$ such that the following inequality holds:
			$$
			\|\nabla g(\phi)\|^2 \leq C(M_{0})\left(1+\|\nabla \Delta \phi\|^{2 \sigma}\right).
			$$
			(ii) If $\|\phi\|_{H^1} \leq M_{0}$,  $|g^{\prime}(x)|<C(|x|^{p}+1)$ $\big(p > 0 \text { arbitrary if } n = 1, 2; $ $ 0 < p < 4 \text { if } n = 3\big)$, and $\left|g^{\prime \prime}(x)\right|<C\left(|x|^{p^{\prime}}+1\right)$ $\big( p^{\prime}>0 \text { arbitrary if } n=1,2;$ $ 0<p^{\prime}<3 \text { if } n=3\big)$, for any $\phi \in H^4$, there exist $0 \leq \sigma<$ 1 and a constant $C(M_{0})$ such that the following inequality holds:
			$$
			\|\Delta g(\phi)\|^2 \leq C(M_{0})\left(1+\left\|\Delta^2 \phi\right\|^{2 \sigma}\right) .
			$$
		\end{lemma}
        \begin{remark}\label{sigma}
        For Allen-Cahn or Cahn-Hilliard equation, $g(x)=k f(x)-s x$ with $k\le M_1$, so that Lemma~\ref{lm_shen} applies with $p=2$ and $p'=1$. Consequently, we have
        \[
        \|\nabla g(\phi)\|^2 \le C(M_0,M_1,s)\left(1+\|\nabla\Delta\phi\|^{2\sigma}\right),
        \]
        where $\sigma$ depends on the spatial dimension $n$: $\sigma=0$ for $n=1$, $\sigma>0$ for $n=2$, and $\sigma=\tfrac12$ for $n=3$. Moreover, it holds that
        \[
        \|\Delta g(\phi)\|^2 \le C(M_0,M_1,s)\left(1+\|\Delta^2\phi\|^{2\sigma}\right),
        \]
        where $\sigma=\tfrac13$ for $n=1$, $\sigma>\tfrac13$ for $n=2$, and $\sigma=\tfrac23$ for $n=3$ (see the proof of Lemma 2.3 in \cite{shen2018convergence}).
        \end{remark}

		{
			\section{Proof of Theorem \ref{th5}}\label{AP2} 
			To prove Theorem \ref{th5}, assuming $r^{n-1}, r^{n}\in [0,M_{1}]$ for $n\le m-1$, we first establish $\|\phi^{m}\|_{H^2}\le M$, where the constant $M$ depends on $T$, $s$, $\phi^{0}$, $\phi^{1}$, $M_{0}$, $M_{1}$, and $\Omega$, and then demonstrate that $r^{m}\in [0,M_{1}]$ when $\alpha\le\gamma_0{\color{black}\tau}$, where $\gamma_0$ depends on $T$, $s$, $\phi^{0}$, $\phi^{1}$, $M_{0}$, $M_{1}$, and $\Omega$.\par
			By multiplying \eqref{AC_nu1} with $4{\color{black}\tau}\phi^{n+1}$ and using integration by parts, then we have
			\begin{align}
				&\|\phi^{n+1}\|^2-\|\phi^{n}\|^2+\|2\phi^{n+1}-\phi^{n}\|^2-\|2\phi^{n}-\phi^{n-1}\|^2+\|\phi^{n+1}-2\phi^{n}+\phi^{n-1}\|^2 \nonumber\\
				&\leq-4{\color{black}\tau}\int_{\Omega}|\nabla\phi^{n+1}|^{2}d\boldsymbol{x}+4{\color{black}\tau}\int_{\Omega}|(\overline{r}^{n+\frac{1}{2}}f(\overline{\phi}^{n+\frac{1}{2}})-s\overline{\phi}^{n+\frac{1}{2}})\phi^{n+1}|d\boldsymbol{x}\nonumber\\
				&-4{\color{black}\tau} s\int_{\Omega}|\phi^{n+1}|^{2}d\boldsymbol{x}.\label{AC0_r_1}
			\end{align}
			{\color{black} By applying the Young’s inequality to the second term on the right-hand side of equation \eqref{AC0_r_1}, we derive}
			\begin{align}
				4{\color{black}\tau}\int_{\Omega}|(\overline{r}^{n+\frac{1}{2}}f(\overline{\phi}^{n+\frac{1}{2}})-s\overline{\phi}^{n+\frac{1}{2}})\phi^{n+1}|d\boldsymbol{x} &\le 4{\color{black}\tau} \frac{1}{s}\int_{\Omega}|(\overline{r}^{n+\frac{1}{2}}f(\overline{\phi}^{n+\frac{1}{2}})-s\overline{\phi}^{n+\frac{1}{2}})|^{2}d\boldsymbol{x}\nonumber\\
				&+{\color{black}\tau} s\int_{\Omega}|\phi^{n+1}|^{2}d\boldsymbol{x}.
			\end{align}
			{\color{black} Since $f'(x) \le M_0$ for all $x\in\mathbb{R}$, $r^{n-1}, r^n \in [0,M_1]$, and
            $\rvert\overline{r}^{n+\frac12} f'(x) - s \rvert \le M_2$, where $M_2$ depends on $M_0$, $M_1$, and $s$. We define $g(\phi) = r^{n+\frac12} f(\phi) - s\phi$. Then, we obtain} 
			\begin{align}
				&4{\color{black}\tau} \frac{1}{s}\int_{\Omega}|g(\overline{\phi}^{n+\frac{1}{2}})|^{2}d\boldsymbol{x}\le 8{\color{black}\tau} \frac{1}{s}\int_{\Omega}|g(\overline{\phi}^{n+\frac{1}{2}})-g(0)|^{2}+|g(0)|^{2}d\boldsymbol{x}\nonumber\\
				&=\frac{8{\color{black}\tau}}{s}\int_{\Omega}|g^{\prime}(\xi)\overline{\phi}^{n+\frac{1}{2}}|^{2}+|g(0)|^{2}d\boldsymbol{x}
				\le 50{\color{black}\tau} \frac{1}{s}M^{2}_{2}\int_{\Omega}|\phi^{n}|^{2}+|\phi^{n-1}|^{2}d\boldsymbol{x}+8{\color{black}\tau} \frac{1}{s}M_{1}|\Omega|.
			\end{align}
			To simplify the notation, we set $C=\max(50\frac{1}{s}M^{2}_{2},8\frac{1}{s}M_{1}|\Omega|)$. Then, we can get
			\begin{align}
				&\|\phi^{n+1}\|^2-\|\phi^{n}\|^2+\|2\phi^{n+1}-\phi^{n}\|^2-\|2\phi^{n}-\phi^{n-1}\|^2+\|\phi^{n+1}-2\phi^{n}+\phi^{n-1}\|^2 \nonumber\\
				&=C{\color{black}\tau}\|\phi^{n}\|^{2}+C{\color{black}\tau}\|\phi^{n-1}\|^{2}+C{\color{black}\tau}.\label{X1X2_L2}
			\end{align}
			Letting $X^{n+1}=\|\phi^{n+1}\|^2+\|2\phi^{n+1}-\phi^{n}\|^2$, we obtain
			$
			X^{n+1}-X^{n}\le {\color{black}\tau} C\left(X^{n}+X^{n-1}\right)+C{\color{black}\tau}.
			$
			By adding ${\color{black}\tau} C X^{n}$ to both sides of the equation, we can get
			\begin{align}
				X^{n+1}+{\color{black}\tau} C X^{n}&\le X^{n} + {\color{black}\tau} C\left(2X^{n}+X^{n-1}\right)+C{\color{black}\tau}\nonumber\\
				&\le \left(1+2{\color{black}\tau} C\right)\left(X^{n}+{\color{black}\tau} C X^{n-1}\right)+C{\color{black}\tau}.
			\end{align}
			By taking the sum from $1$ to $m-1$,
			\begin{align}
				\left(X^{m}+{\color{black}\tau} C X^{m-1}\right)\le \left(X^{2}+{\color{black}\tau} C X^{1}\right)+CT+2{\color{black}\tau} C\sum_{n=1}^{m-1}\left(X^{n}+{\color{black}\tau} C X^{n-1}\right).
			\end{align}
			By Lemma \ref{lm3}, we can get
			\begin{align}
				\left(X^{m}+{\color{black}\tau} C X^{m-1}\right)\le \left(X^{2}+{\color{black}\tau} C X^{1}+CT\right)e^{2TC}.
			\end{align}
			{Using the first-order scheme \eqref{AC_1st_1}--\eqref{AC_1st_2} and following a similar proof process as above, along with using \eqref{X1X2_L2}, we can easily establish the boundedness of $X^{1}$ and $X^{2}$.} Then, we have $\|\phi^{m}\|^{2}_{L^{2}}\le \overline{M_{0}},~\forall m\geq 2$, where $\overline{M_{0}}$ depends on $M_{0}$, $M_{1}$, $s$, $T$, $\Omega$, $\phi^0$, and $\phi^1$.\par
			Next, we will prove $\|\Delta\phi^{m}\|\le \overline{M}_{1}$, where $\overline{M}_{1}$ depends on $M_{0}$, $M_{1}$, $s$, $T$, $\Omega$, $\phi^0$, and $\phi^1$. By multiplying \eqref{AC_nu1} with $4{\color{black}\tau} \Delta^{2}\phi^{n+1}$ and using integration by parts, then we have
			\begin{align}
				&\|\Delta\phi^{n+1}\|^2-\|\Delta\phi^{n}\|^2+\|2\Delta\phi^{n+1}-\Delta\phi^{n}\|^2-\|2\Delta\phi^{n}-\Delta\phi^{n-1}\|^2+\|\Delta\phi^{n+1}-2\Delta\phi^{n}\nonumber\\
				&+\Delta\phi^{n-1}\|^2 =-4{\color{black}\tau}\int_{\Omega}|\nabla\Delta\phi^{n+1}|^{2}d\boldsymbol{x}-4{\color{black}\tau} s\int_{\Omega}|\Delta\phi^{n+1}|^{2}d\boldsymbol{x}+\nonumber\\
				&4{\color{black}\tau}\int_{\Omega}\nabla(\overline{r}^{n+\frac{1}{2}}f(\overline{\phi}^{n+\frac{1}{2}})-s\overline{\phi}^{n+\frac{1}{2}})\cdot\nabla\Delta\phi^{n+1}d\boldsymbol{x} .\label{AC_r_2}
			\end{align}
			By applying the Young's inequality to the final term on the right-hand side of equation \eqref{AC_r_2}, we derive
			\begin{align}
				&4{\color{black}\tau}\int_{\Omega}\nabla(\overline{r}^{n+\frac{1}{2}}f(\overline{\phi}^{n+\frac{1}{2}})-s\overline{\phi}^{n+\frac{1}{2}})\cdot\nabla\Delta\phi^{n+1}d\boldsymbol{x} \nonumber\\
				&\le 4{\color{black}\tau}\int_{\Omega}|\nabla(\overline{r}^{n+\frac{1}{2}}f(\overline{\phi}^{n+\frac{1}{2}})-s\overline{\phi}^{n+\frac{1}{2}})|^{2}d\boldsymbol{x}+{\color{black}\tau}\int_{\Omega}|\nabla\Delta\phi^{n+1}|^{2}d\boldsymbol{x}.
			\end{align}
			Because $|\overline{r}^{n+\frac{1}{2}}f^{\prime}(x)-s|<M_{2}$ and by using the Gagliardo–Nirenberg inequality, we obtain 
			\begin{align}
				&4{\color{black}\tau} \int_{\Omega}|\nabla(\overline{r}^{n+\frac{1}{2}}f(\overline{\phi}^{n+\frac{1}{2}})-s\overline{\phi}^{n+\frac{1}{2}})|^{2}d\boldsymbol{x} \le 4{\color{black}\tau} M^{2}_{2}\int_{\Omega}|\nabla\overline{\phi}^{n+\frac{1}{2}}|^{2}d\boldsymbol{x}\nonumber\\
				&\le 4{\color{black}\tau} M^{2}_{2}C_{1}\left(\int_{\Omega}|\Delta\overline{\phi}^{n+\frac{1}{2}}|^{2}d\boldsymbol{x}+\|\overline{\phi}^{n+\frac{1}{2}}\|^{2}\right)\nonumber\\
				&\le 30{\color{black}\tau} M^{2}_{2}C_{1}\left(\int_{\Omega}|\Delta\phi^{n}|^{2}d\boldsymbol{x} + \int_{\Omega}|\Delta\phi^{n-1}|^{2}d\boldsymbol{x}+\overline{M_{0}}\right),
			\end{align}
			where the constant $C_1\ge 0$ depends on the domain $\Omega$. To simplify the notation, we set $\overline{C}=30M^{2}_{2}C_{1}$. Then, we can get
			\begin{align}
				&\|\Delta\phi^{n+1}\|^2-\|\Delta\phi^{n}\|^2+\|2\Delta\phi^{n+1}-\Delta\phi^{n}\|^2-\|2\Delta\phi^{n}-\Delta\phi^{n-1}\|^2 \nonumber\\
				&\le{\color{black}\tau} \overline{C}\left(\|\Delta\phi^{n}\|^{2} + \|\Delta\phi^{n-1}\|^{2}\right)+{\color{black}\tau}\overline{C}\overline{M_{0}}.
			\end{align}
			Letting $X^{n+1}=\|\Delta\phi^{n+1}\|^2+\|2\Delta\phi^{n+1}-\Delta\phi^{n}\|^2$ and 
			adding ${\color{black}\tau} \overline{C} X^{n}$ to both sides of the equation, we obtain
			\begin{align}
				&X^{n+1}+{\color{black}\tau} \overline{C} X^{n}\le X^{n} + {\color{black}\tau} \overline{C}\left(2X^{n}+X^{n-1}\right)+{\color{black}\tau}\overline{C}\overline{M_{0}}\nonumber\\
				&\le\left(1+2{\color{black}\tau} \overline{C}\right)\left(X^{n}+{\color{black}\tau} \overline{C} X^{n-1}\right)+{\color{black}\tau}\overline{C}\overline{M_{0}}.
			\end{align}
			By taking the sum from $1$ to $m-1$,
			\begin{align}
				\left(X^{m}+{\color{black}\tau} \overline{C} X^{m-1}\right)\le \left(X^{2}+{\color{black}\tau} \overline{C} X^{1}+T\overline{C}\overline{M_{0}}\right)+2{\color{black}\tau}\overline{C}\sum_{n=1}^{m-1}\left(X^{n}+{\color{black}\tau} \overline{C} X^{n-1}\right).
			\end{align}
			By Lemma \ref{lm3}, we can get
			\begin{align}
				\left(X^{m}+{\color{black}\tau} \overline{C} X^{m-1}\right)\le \left(X^{2}+{\color{black}\tau} \overline{C} X^{1}+T\overline{C}\overline{M_{0}}\right)e^{2T\overline{C}}.
			\end{align}
			{Similar to the above text, we can easily establish the boundedness of $X^{1}$ and $X^{2}$.} Then, we can get $\|\Delta\phi^{m}\|\le \overline{M}_{1}$, where $\overline{M}_{1}$ depends on $M_{0}$, $M_{1}$, $s$, $T$, $\Omega$, $\phi^0$, and $\phi^1$. By using Gagliardo–Nirenberg interpolation inequality, we can get $\|\nabla\phi^{m}\|\le \tilde{C}_{0}\overline{M}_{1}$. It follows that $\|\phi^{m}\|_{H^{2}}\le M,~\forall m\geq 2$, where $M$ depends on $M_{0}$, $M_{1}$, $s$, $T$, $\Omega$, $\phi^0$, and $\phi^1$.
			
			Now, we can get $\|\phi^{n+1}\|_{H^2}\le M$, $\forall n\le m-1$. Since the dimension $d\leq 3$, it is easy to get $\|\phi^{n+1}\|_{L^{\infty}}\le C_\Omega M$ with $C_\Omega$ the Sobolev embedding constant. Therefore, $\exists M_{3}$, s.t. $E_{0}(\phi^{n+1})$, $E_{0}(\phi^{n})$, $E_{0}(\phi^{n-1})$, $f(\overline{\phi}^{n+\frac{1}{2}})\le M_{3}$, where $M_3$ depends on $C_\Omega$ and $M$. By assumption, $r^{n-1}, r^{n}\in [0,M_{1}]$ for $n\le m-1$. {Let $z^{n}=|\frac{3r^{n}-r^{n-1}}{2}-1|$, $1\le n\le m-1$.} From \eqref{AC_nu2}, we can get 
			\begin{align}
				&z^{n+1}-z^{n}\le \left|\left(\frac{3r^{n+1}-r^{n}}{2}-1\right)-\left(\frac{3r^{n}-r^{n-1}}{2}-1\right)\right|\nonumber\\ &=\alpha\bigg|-\frac{3E^{n+1}_{0}-4E^{n}_{0}+E^{n-1}_{0}}{2}
				+\int_{\Omega}(\overline{r}^{n+\frac{1}{2}}f(\overline{\phi}^{n+\frac{1}{2}})-s\overline{\phi}^{n+\frac{1}{2}})\frac{3\phi^{n+1}-4\phi^{n}+\phi^{n-1}}{2}d\boldsymbol{x}\bigg|\nonumber\\
				&\le\alpha\bigg|-\frac{3E^{n+1}_{0}-4E^{n}_{0}+E^{n-1}_{0}}{2}+\int_{\Omega}({(2-r^{n})}f(\overline{\phi}^{n+\frac{1}{2}})-s\overline{\phi}^{n+\frac{1}{2}})\frac{3\phi^{n+1}-4\phi^{n}+\phi^{n-1}}{2}d\boldsymbol{x}\bigg|\nonumber\\
				&\quad+\alpha{2z^{n}}\int_{\Omega}\big|f(\overline{\phi}^{n+\frac{1}{2}})\frac{3\phi^{n+1}-4\phi^{n}+\phi^{n-1}}{2}\big|d\boldsymbol{x}\nonumber\\
				&\le\alpha C(M,M_{1},M_{2},M_{3},\Omega,s)+\alpha z^{n}C(M,M_{1},M_{2},M_{3},\Omega,s). \label{zn}
			\end{align}
			Taking the summation of \eqref{zn} for $n$ from $1$ to $m-1$,
			\begin{align}
				z^{m}-z^{1}\le \sum_{n=1}^{m-1}C(M,M_{1},M_{2},M_{3},\Omega,s)\alpha + \sum_{n=1}^{m-1}\alpha C(M,M_{1},M_{2},M_{3},\Omega,s) z^{n}.
			\end{align}
			By setting $\alpha={\color{black}\tau}\gamma$, we can get
			\begin{align}
				z^{m}\le z^{1} + TC(M,M_{1},M_{2},M_{3},\Omega,s)\gamma + {\color{black}\tau}\sum_{n=1}^{m-1}\gamma C(M,M_{1},M_{2},M_{3},\Omega,s) z^{n}.
			\end{align}
			By applying Lemma \ref{lm3}, we can get
			\begin{align}
				z^{m}\le e^{T\overline{C}\gamma}\left(z^{1} + T\overline{C}\gamma\right),
			\end{align}
			where $\overline{C}= C(M,M_{1},M_{2},M_{3},\Omega,s)$. {Using the first-order scheme \eqref{AC_1st_1}--\eqref{AC_1st_2} and following a similar proof process as above, we can easily establish the boundedness of $z^{1}$.} Then, it follows that $ \lim_{\gamma \to 0}z^{m}=0$. Therefore, there exists a $\gamma_{0}$ such that $z^{m}<\min(1,M_{1}-1)$ for any $\gamma\leq \gamma_0$, which indicates $r^{m}\in [0,M_{1}]$.
			
			\section{Proof of Theorem \ref{th6}}\label{AP3}
			  To prove Theorem \ref{th6}, assuming $r^{n-1}, r^{n}\in [0,M_{1}]$ for $n\le m-1$, we first establish $\|\phi^{m}\|_{H^2}\le M$, where the constant $M$ depends on $T$, $s$, $\phi^{0}$, $\phi^{1}$, $M_{0}$, $M_{1}$, and $\Omega$, and then demonstrate that $r^{m}\in [0,M_{1}]$ when $\alpha\le\gamma_0{\color{black}\tau}$, where $\gamma_0$ depends on $T$, $s$, $\phi^{0}$, $\phi^{1}$, $M_{0}$, $M_{1}$, and $\Omega$.\par
			We multiply \eqref{CH_nu1} with $-4{\color{black}\tau} \Delta\phi^{n+1}$ and use integration by parts, then we can get
			\begin{align}
				&\|\nabla\phi^{n+1}\|^2-\|\nabla\phi^{n}\|^2+\|2\nabla\phi^{n+1}-\nabla\phi^{n}\|^2-\|2\nabla\phi^{n}-\nabla\phi^{n-1}\|^2\nonumber\\
				&+\|\nabla\phi^{n+1}-2\nabla\phi^{n}+\nabla\phi^{n-1}\|^2 =-4{\color{black}\tau}\int_{\Omega}|\nabla\Delta\phi^{n+1}|^{2}d\boldsymbol{x}-4{\color{black}\tau} s\int_{\Omega}|\Delta\phi^{n+1}|^{2}d\boldsymbol{x}\nonumber\\
				&+4{\color{black}\tau}\int_{\Omega}\nabla(\overline{r}^{n+\frac{1}{2}}f(\overline{\phi}^{n+\frac{1}{2}})-s\overline{\phi}^{n+\frac{1}{2}})\cdot\nabla\Delta\phi^{n+1}d\boldsymbol{x} .\label{AC_r_1}
			\end{align}
			Similar to the arguments in the last Appendix, we can get
			\begin{align}
				&\|\nabla\phi^{n+1}\|^2-\|\nabla\phi^{n}\|^2+\|2\nabla\phi^{n+1}-\nabla\phi^{n}\|^2-\|2\nabla\phi^{n}-\nabla\phi^{n-1}\|^2\nonumber\\
				&\le {\color{black}\tau} C\big(\|\nabla\phi^{n}\|^{2}+ \|\nabla\phi^{n-1}\|^{2}\big),
			\end{align}
			where $C=20C_{0}M^{2}_{2}$. Let $X^{n+1}=\|\nabla\phi^{n+1}\|^2+\|2\nabla\phi^{n+1}-\nabla\phi^{n}\|^2$ and add ${\color{black}\tau} C X^{n}$ to both sides of the equation. Then, we obtain
			\begin{align}
				X^{n+1}+{\color{black}\tau} C X^{n}\le X^{n} + {\color{black}\tau} C\left(2X^{n}+X^{n-1}\right)\le \left(1+2{\color{black}\tau} C\right)\left(X^{n}+{\color{black}\tau} C X^{n-1}\right).
			\end{align}
			By taking the sum from $1$ to $m-1$,
			\begin{align}
				\left(X^{m}+{\color{black}\tau} C X^{m-1}\right)\le \left(X^{2}+{\color{black}\tau} C X^{1}\right)+2{\color{black}\tau} C\sum_{n=1}^{m-1}\left(X^{n}+{\color{black}\tau} C X^{n-1}\right).
			\end{align}
            {\color{black} By Lemma~4, we obtain
            \begin{align}
            (X^m + \tau C X^{m-1}) \le (X^2 + \tau C X^1)e^{2TC}.
            \end{align}
             It follows that $\|\nabla \phi^m\| \le \overline{M_1}$, where the constant $\overline{M_1}$ depends on $M_0$, $M_1$, $s$, $\phi^0$, $\phi^1$, $\Omega$, and $T$. Using the Poincaré inequality together with the mass conservation property $\int_\Omega \phi^m \, dx = \int_\Omega \phi^0 \, dx$, we have}
			\begin{align}
				{\color{black}\|\phi^{m}\|\le C\|\nabla\phi^{m}\|}+\|\frac{1}{|\Omega|}\int_{\Omega}\phi^{m}d\boldsymbol{x}\|\le C\overline{M}_{1}+\int_{\Omega}\phi^{0}d\boldsymbol{x}.
			\end{align}
			Therefore, there exists a constant $\overline{M}_2$, depending on $\phi^0$ and $\overline{M}_1$, such that $\|\phi^{n+1}\|_{H^1} \le \overline{M}_2$, $\forall n \le m-1$. \par
			Next we will proof $\|\Delta \phi^{m}\|\le \overline{M}_{3}$, where $\overline{M}_{3}$ depends on $T$, $s$, $\phi^{0}$, $\phi^{1}$, $M_{0}$, $M_{1}$, and $|\Omega|$. {\color{black} By multiplying \eqref{CH_nu1} with $4{\color{black}\tau} \Delta^{2}\phi^{n+1}$, we have}
			\begin{align}
				&\|\Delta\phi^{n+1}\|^2-\|\Delta\phi^{n}\|^2+\|2\Delta\phi^{n+1}-\Delta\phi^{n}\|^2-\|2\Delta\phi^{n}-\Delta\phi^{n-1}\|^2+\|\Delta\phi^{n+1}-2\Delta\phi^{n}\nonumber\\
				&+\Delta\phi^{n-1}\|^2+4{\color{black}\tau}\int_{\Omega}|\Delta^{2}\phi^{n+1}|^{2}d\boldsymbol{x}+4{\color{black}\tau} s\int_{\Omega}|\nabla\Delta\phi^{n+1}|^{2}d\boldsymbol{x}=\nonumber\\ 
				&4{\color{black}\tau}\int_{\Omega}\Delta(\overline{r}^{n+\frac{1}{2}}f(\overline{\phi}^{n+\frac{1}{2}})-s\overline{\phi}^{n+\frac{1}{2}})\Delta^{2}\phi^{n+1}d\boldsymbol{x}\leq 4{\color{black}\tau} \overline{C}\left\|\Delta (\overline{r}^{n+\frac{1}{2}}f(\overline{\phi}^{n+\frac{1}{2}})-s\overline{\phi}^{n+\frac{1}{2}})\right\|^2\nonumber\\
				&+{\color{black}\tau}\left\|\Delta^2 \phi^{n+1}\right\|^2.\label{CH_r_1}
			\end{align}
			By the Young's inequality and Lemma \ref{lm_shen}, for any $\zeta>0$, there exists a constant \\{\color{black} $C_1(M_1,s,\zeta, \sigma, \overline{M}_{2})$} depending on $\zeta$, such that the following inequality holds:
			\begin{align}
                &\overline{C}\left\|\Delta (\overline{r}^{n+\frac{1}{2}}f(\overline{\phi}^{n+\frac{1}{2}})-s\overline{\phi}^{n+\frac{1}{2}})\right\|^2 \leq C_0(M_1,s,\overline{M}_{2})\left(1+\left\|\Delta^2 \overline{\phi}^{n+\frac{1}{2}}\right\|^{2 \sigma}\right) \nonumber\\
                &= C_0(M_1,s,\overline{M}_{2})\left\|\Delta^2 \overline{\phi}^{n+\frac{1}{2}}\right\|^{2\sigma}+C_0(M_1,s,\overline{M}_{2})\nonumber\\
                &\leq \zeta\sigma\left\|\Delta^2 \overline{\phi}^{n+\frac{1}{2}}\right\|^2+{\color{black} \zeta^{-\frac{\sigma}{1-\sigma}}(1-\sigma)C_0(M_1,s,\overline{M}_{2})^{\frac{1}{1-\sigma}}+C_0(M_1,s,\overline{M}_{2})}\nonumber\\
                &\leq 8\sigma\zeta\left(\left\|\Delta^2 \phi^{n}\right\|^2+\left\|\Delta^2 \phi^{n-1}\right\|^2\right)+{\color{black} C_1(M_1,s,\zeta, \sigma, \overline{M}_{2})},
            \end{align}
			{\color{black}where the value of $\sigma$ is specified in Remark~\ref{sigma}.} We choose $\zeta = 1/(32\sigma)$ to obtain
			\begin{align}
				&\|\Delta\phi^{n+1}\|^2-\|\Delta\phi^{n}\|^2+\|2\Delta\phi^{n+1}-\Delta\phi^{n}\|^2-\|2\Delta\phi^{n}-\Delta\phi^{n-1}\|^2+\|\Delta\phi^{n+1}-2\Delta\phi^{n}\nonumber\\
				&+\Delta\phi^{n-1}\|^2 +3{\color{black}\tau}\|\Delta^{2}\phi^{n+1}\|^{2}-{\color{black}\tau}\left(\left\|\Delta^2 \phi^{n}\right\|^2+\left\|\Delta^2 \phi^{n-1}\right\|^2\right)\leq {\color{black}\tau} C_2(M_1,s,\sigma,\overline{M}_{2}).
			\end{align}
			Taking the sum from 0 to m-1, we can get
			\begin{align}
				&\|\Delta\phi^{m}\|^2-\|\Delta\phi^{1}\|^2+\|2\Delta\phi^{m}-\Delta\phi^{m-1}\|^2-\|2\Delta\phi^{1}-\Delta\phi^{0}\|^2 \nonumber\\
				&+\sum_{n=2}^{m}{\color{black}\tau}\|\Delta^{2}\phi^{n}\|^{2}-{\color{black}\tau}\left(\left\|\Delta^2 \phi^{1}\right\|^2+\left\|\Delta^2 \phi^{0}\right\|^2\right)\leq TC_2(M_1,s,\sigma,\overline{M}_{2}).
			\end{align}
			{Using the first-order scheme \eqref{CH_1st_1}--\eqref{CH_1st_2}, we can easily establish the boundedness of $\left\|\Delta^2 \phi^{1}\right\|^2$.} Then, we can get the conclusion.
		}
		\section{Proof of Theorem \ref{th7}}\label{AP4}
		To prove Theorem \ref{th7}, we first state and prove a $H^3$-norm stability result of the numerical solution. 
		\begin{theorem}[$H^3$-norm stability of $\phi^n$]\label{th9}
			Assume the initial value $\phi^{0}\in H^{3}(\Omega)$. Then for all $n\le T/{\color{black}\tau}$, the solution $\phi^{n}$ in \eqref{AC_1st_1}--\eqref{AC_1st_2} satisfies
			\begin{align}
				&\|\Delta\phi^{n+1}\|^{2}+{\color{black}\tau}\|\nabla\Delta\phi^{n+1}\|^{2}\le C(T,M) +\|\Delta\phi^{0}\|^{2}+{\color{black}\tau}\|\nabla\Delta\phi^{0}\|^{2}.
			\end{align}
		\end{theorem}
		\begin{proof} Based on Theorems \ref{AP1th2}, through induction, we establish that $\|\phi^{n}\|_{\infty}\le1$ and $r^{n}\le M_1$ for all $n$, where $M_1>1$. Because $r_n > 0$, $\forall 0 \leq n \leq T / {\color{black}\tau}$, we can conclude by the energy stability in  Theorem \ref{th2} that $\phi^{n} \in H^1(\Omega)$, $\exists M>0, \|\phi^n\|_{H^1}\le M$, $\forall 0 \leq n \leq T / {\color{black}\tau}$.
			
			By multiplying \eqref{AC_1st_1} with $\Delta^{2}\phi^{n+1}$, we can get
			\begin{align}
				\frac{\Delta^{2}\phi^{n+1}(\phi^{n+1}-\phi^{n})}{{\color{black}\tau}}=\Delta\phi^{n+1}\Delta^{2}\phi^{n+1}-s\phi^{n+1}\Delta^{2}\phi^{n+1}-(r^{n}f(\phi^{n})-s\phi^{n})\Delta^{2}\phi^{n+1}.
			\end{align}
			Then we can get
			\begin{align}
				&\frac{\|\Delta\phi^{n+1}\|^{2}-\|\Delta\phi^{n}\|^{2}+\|\Delta\phi^{n+1}-\Delta\phi^{n}\|^{2}}{2{\color{black}\tau}}=-\|\nabla\Delta\phi^{n+1}\|^{2}-s\|\Delta\phi^{n+1}\|^{2}+\nonumber\\
				&\big((r^{n}\nabla f(\phi^{n})-s\nabla\phi^{n}),\nabla\Delta\phi^{n+1}\big)\le-\|\nabla\Delta\phi^{n+1}\|^{2}-s\|\Delta\phi^{n+1}\|^{2}\nonumber\\
				&+ C\|(r^{n}\nabla f(\phi^{n})-s\nabla\phi^{n})\|^{2}+\frac{1}{2}\|\nabla\Delta\phi^{n+1}\|^{2}.
			\end{align}
			{\color{black} By Lemma \ref{lm4} and the Young's inequality,
            \begin{align}
                &C\left\|\left(r^n \nabla f\left(\phi^n\right)-s \nabla \phi^n\right)\right\|^2 \leq C_0(M_1,s,M)\left(1+\left\|\nabla\Delta \phi^{n}\right\|^{2 \sigma}\right) \nonumber\\
                &= C_0(M_1,s,M)\left\|\nabla\Delta \phi^{n}\right\|^{2\sigma}+C_0(M_1,s,M)\nonumber\\
                &\leq \zeta\sigma\left\|\nabla\Delta \phi^{n}\right\|^2+{ \zeta^{-\frac{\sigma}{1-\sigma}}(1-\sigma)C_0(M_1,s,M)^{\frac{1}{1-\sigma}}+C_0(M_1,s,M)},
            \end{align}
            where the value of $\sigma$ is specified in Remark~\ref{sigma}. If $\sigma = 0$, it is trivial to see that the above inequality holds. Otherwise, by taking $\zeta = 1/(2\sigma)$, we obtain $C\|(r^{n}\nabla f(\phi^{n})-s\nabla\phi^{n})\|^{2}\le \frac{1}{2}\|\nabla\Delta\phi^{n}\|^{2}+C_1(M_1,s,\sigma,M)$, which yields}
			\begin{align}
				&\|\Delta\phi^{n+1}\|^{2}-\|\Delta\phi^{n}\|^{2}+{\color{black}\tau}\|\nabla\Delta\phi^{n+1}\|^{2}-{\color{black}\tau}\|\nabla\Delta\phi^{n}\|^{2}\le 2{\color{black}\tau} C_1(M_1,s,\sigma,M).
			\end{align}
			We conclude the proof by taking the sum from $0$ to $n-1$.
		\end{proof}
		
		We are now ready to prove Theorem \ref{th7}.
		By setting $w^{n+1}=\mu^{n+1}-\mu(t^{n+1})$, $g(\phi^{n})=f(\phi^{n})-s\phi^{n}$. 
		Because $\phi(t_{0})\in H^{2}(\Omega)$, we know $\phi(t)\in C([0,T];H^{2}(\Omega))$ (see \cite[Theorem 2.6]{shen2018convergence}). By Theorems \ref{th2}, \ref{AP1th2}, and \ref{th9}, we can deduce that $\exists M_2>0, \|\phi^n\|_{H^2}\le M_2$, $\forall 0 \leq n \leq T / {\color{black}\tau}$, and $\phi^{n}\in H^{2}(\Omega)\subsetneqq L^{\infty}(\Omega)$, we can find a constant $\overline{C}$ such that $|g(\phi)|$, $|g^{\prime}(\phi)|$, $|g(\phi^{n})|$, $|f(\phi^{n})|$, $|g^{\prime}(\phi^{n})|$, $|E_{0}(\phi)|$, $|E_{0}(\phi^{n})|\le \overline{C}$.\par
        {\color{black} We define $z^{n+1} = |r^{n+1} - 1|$. Using \eqref{AC_1st_2}, we obtain}
		\begin{align}
			&z^{n+1}-z^{n}=(|r^{n+1}-1|-|r^{n}-1|)\le|(r^{n+1}-1)-(r^{n}-1)| \nonumber\\
			&=\alpha\bigg|-\big(E_{0}(\phi^{n+1})-E_{0}(\phi^{n})\big)+(r^{n}-1)\int_{\Omega}f(\phi^{n})\big(\phi^{n+1}-\phi^{n}\big)d\boldsymbol{x}\nonumber\\
			&+\int_{\Omega}\left(f(\phi^{n})-s\phi^{n}\right)\big(\phi^{n+1}-\phi^{n}\big)d\boldsymbol{x}\bigg|\nonumber\\
			&\le\alpha\bigg|\big(E_{0}(\phi^{n+1})-E_{0}(\phi^{n})\big)-\int_{\Omega}g(\phi^{n})\big(\phi^{n+1}-\phi^{n}\big)d\boldsymbol{x}\bigg|\nonumber\\
            &+\alpha z^{n}\int_{\Omega}\left|f(\phi^{n})\big(\phi^{n+1}-\phi^{n}\big)\right|d\boldsymbol{x}.  \label{zn_1}
		\end{align}
		{\color{black}For the first term of the right hand side of \eqref{zn_1}, since} $E_{0}(\phi)=\int_{\Omega}(F(\phi)-\frac{s}{2}\phi^{2})d\boldsymbol{x}$ and $\frac{d}{d\phi}(F(\phi)-\frac{s}{2}\phi^{2})=g(\phi)$, we can obtain 
		\begin{align}
			&\bigg|\big(E_{0}(\phi^{n+1})-E_{0}(\phi^{n})\big)-\int_{\Omega}g(\phi^{n})\big(\phi^{n+1}-\phi^{n}\big)d\boldsymbol{x}\bigg|\nonumber \\
			&=\bigg|\int_{\Omega}\left(F(\phi^{n+1})-\frac{s}{2}(\phi^{n+1})^{2}\right)-\left(F(\phi^{n})-\frac{s}{2}(\phi^{n})^{2}\right)-g(\phi^{n})\big(\phi^{n+1}-\phi^{n}\big)d\boldsymbol{x}\bigg|\nonumber \\
			&=\bigg|\int_{\Omega}g(\xi_{1})\big(\phi^{n+1}-\phi^{n}\big)-g(\phi^{n})\big(\phi^{n+1}-\phi^{n}\big)d\boldsymbol{x}\bigg|\nonumber \\  
			&=\bigg|\int_{\Omega}g^{\prime}(\xi_{2})\big(\xi_{1}-\phi^{n}\big)\big(\phi^{n+1}-\phi^{n}\big)d\boldsymbol{x}\bigg|\le \overline{C}\int_{\Omega}\big|\xi_{1}-\phi^{n}\big|\big|\phi^{n+1}-\phi^{n}\big|d\boldsymbol{x} .\label{int_E_g}
		\end{align}
		{Because $\xi_{1}$ is between $\phi^{n+1}$ and $\phi^{n}$, and $\phi^{n} \in H^{2}(\Omega)$, we can get
			\begin{align}
				&\int_{\Omega}\big|\xi_{1}-\phi^{n}\big|\big|\phi^{n+1}-\phi^{n}\big|d\boldsymbol{x}\le\int_{\Omega}\big|\phi^{n+1}-\phi^{n}\big|^{2}d\boldsymbol{x}\nonumber\\
				&={\color{black}\tau}^{2}\int_{\Omega}\big|\Delta \phi^{n+1}-s\phi^{n+1}-(r^{n}f(\phi^{n})-s\phi^{n})\big|^{2}d\boldsymbol{x}\le C_{1}{\color{black}\tau}^{2}, \label{xi_phi}
			\end{align}
			where $r^{n}\in [0,M_1]$, $\phi^{n} \in H^{2}(\Omega) \hookrightarrow L^\infty(\Omega)$, and the constant $C_1$ depends on $M_1$, $M_2$, $\overline{C}$, and $\Omega$.}
        {\color{black}For the second term on the right-hand side of \eqref{zn_1}, using the Hölder inequality and the fact that $\phi^{n} \in H^{2}(\Omega)\hookrightarrow L^\infty(\Omega)$, we obtain 
        \begin{align}
            &\alpha z^{n}\int_{\Omega}\left|f(\phi^{n})\big(\phi^{n+1}-\phi^{n}\big)\right|d\boldsymbol{x}\leq \alpha z^{n}\overline{C}\int_{\Omega}\big|\phi^{n+1}-\phi^{n}\big|d\boldsymbol{x}\nonumber\\
            &\leq \alpha z^{n}\overline{C}\widetilde{C}\|\phi^{n+1}-\phi^{n}\|^2=\alpha z^{n}\overline{C}\widetilde{C}\tau\|\Delta \phi^{n+1}-s\phi^{n+1}-(r^{n}f(\phi^{n})-s\phi^{n})\|^2\nonumber\\
            &\leq C_2\alpha{\color{black}\tau} z^{n}.
        \end{align}}    
		By combining \eqref{zn_1}, \eqref{int_E_g} and \eqref{xi_phi}, we can get
		\begin{align}
			&z^{n+1}-z^{n}\le C_2\alpha{\color{black}\tau} z^{n}+C_1\alpha {\color{black}\tau}^{2}. \label{zn_2}
		\end{align}
		Taking the summation of \eqref{zn_1} for $n$ from $0$ to $m$, 
		$
		z^{m+1}\le C_2{\color{black}\tau}\alpha \sum_{n=0}^{m}  z^{n}+C_1{\color{black}\tau}\alpha\le C_3{\color{black}\tau}\alpha \sum_{n=0}^{m}  z^{n}+C_3{\color{black}\tau}\alpha,
		$
		where $C_3=max(C_1,C_2)$. Using the discrete Gronwall lemma, we can get
		\begin{align}
			&z^{m}\le C_3\exp((1-{\color{black}\tau}\alpha C_3)^{-1}t^{m}C_3\alpha)\alpha{\color{black}\tau} .\label{zn_4}
		\end{align}
		(i) and (ii) in Theorem \ref{th7} then follow from the above inequality.

\bibliographystyle{elsarticle-num.bst}
\bibliography{refs}
\end{document}